\documentclass[11pt,reqno,a4paper]{amsart}
\usepackage{amsmath,amsthm,amsfonts,amssymb,amscd,amstext}
\usepackage[ansinew]{inputenc}
\usepackage{graphicx}
\usepackage{euscript}
\usepackage[normalem]{ulem} % either use this (simple) or
\usepackage{a4wide}

\usepackage{color}
\usepackage[colorlinks=true,backref]{hyperref}

\numberwithin{equation}{section}

\newcommand{\qand}{\quad\text{and}\quad}

\theoremstyle{plain}
\newtheorem{maintheorem}{Theorem}
\newtheorem{maincorollary}[maintheorem]{Corollary}

\newtheorem{theorem}{Theorem}[section]
\newtheorem{proposition}[theorem]{Proposition}
\newtheorem{corollary}[theorem]{Corollary}
\newtheorem{lemma}[theorem]{Lemma}

\theoremstyle{definition}
\newtheorem{remark}[theorem]{Remark}

\newcommand{\RR}{{\mathbb R}}
\newcommand{\BB}{{\mathbb B}}
\newcommand{\NN}{{\mathbb N}}
\newcommand{\ZZ}{{\mathbb Z}}
\newcommand{\FF}{{\mathbb F}}
\newcommand{\EE}{{\mathbb E}}
\newcommand{\TT}{{\mathbb T}}

\newcommand{\sS}{{\mathbb S}}

\newcommand{\vfi}{\varphi}

\newcommand{\de}{\delta}
\newcommand{\De}{\Delta}

\newcommand{\var}{\operatorname{var}}
\newcommand{\diam}{\operatorname{diam}}
\renewcommand{\epsilon}{\varepsilon}
\newcommand{\dist}{\operatorname{dist}}

\newcommand{\esup}{\operatorname{ess\,sup}}

\newcommand{\Leb}{\operatorname{Leb}}

\newcommand{\wt}{\widetilde}

\newcommand{\supp}{\operatorname{supp}}

\newcommand{\cE}{\EuScript{E}}

\newcommand{\D}{\EuScript{D}}
\newcommand{\cP}{\EuScript{P}}

\newcommand{\U}{\EuScript{U}}

\newcommand{\cC}{\EuScript{C}}
\newcommand{\cS}{\EuScript{S}}
\newcommand{\M}{\EuScript{M}}
\newcommand{\B}{\EuScript{B}}
\newcommand{\R}{\EuScript{R}}

\newcommand{\mC}{\mathcal{C}}

\title[Large deviations for N.U.E. maps]
{Large deviations for non-uniformly expanding maps}

\author{V. Ara\'ujo and M. J. Pacifico}

\date{\today}
\begin{thanks} {V.A. was partially supported by CNPq-Brazil
    and FCT-Portugal through CMUP and POCI/MAT/61237/2004.
    M.J.P. was partially supported by
    CNPq-Brazil/FAPERJ-Brazil/Pronex Dynamical Systems.}
\end{thanks}
\begin{document}

% \address{José F. Alves, Centro de Matemática da
%   Universidade do Porto
% Rua do Campo Alegre 687, 4169-007 Porto, Portugal}
% \email{jfalves@fc.up.pt} %\urladdr{http://www.fc.up.pt/cmup/jfalves}

\address{ V\'\i tor Ara\'ujo, Instituto de Matem\'a\-tica e
  Estat\'{\i}stica, Universidade Federal da Bahia, Av.\
  Ademar de Barros s/n, 40170-110 Salvador, Brazil}
  % Universidade Federal do Rio de Janeiro, C. P. 68.530,
  % 21.945-970 Rio de Janeiro, RJ-Brazil \emph{and} Centro de
  % Matem\'atica da Universidade do Porto, Rua do Campo Alegre
  % 687, 4169-007 Porto, Portugal.}
% Departamento de Matem\'atica, PUC-Rio,
%   Rua Marqu\^es de S. Vicente, 225, G\'avea,
%   22453-900,

\email{vitor.araujo.im.ufba@gmail.com \text{and} vitor.d.araujo@ufba.br}
%\urladdr{http://www.fc.up.pt/cmup/home/vdaraujo}

\address{Maria Jos\'e Pacifico,
Instituto de Matem\'atica,
Universidade Federal do Rio de Janeiro,
C. P. 68.530, 21.945-970 Rio de Janeiro, Brazil}
\email{pacifico@im.ufrj.br}

% \address{Vilton  Pinheiro, Departamento de Matem\'atica, Universidade Federal da Bahia\\
% Av. Ademar de Barros s/n, 40170-110 Salvador, Brazil.}
% \email{viltonj@ufba.br}


\renewcommand{\subjclassname}{\textup{2000} Mathematics Subject Classification}


\begin{abstract}
  We obtain large deviation bounds for non-uniformly
  expanding maps with non-flat singularities or
  criticalities and for partially hyperbolic non-uniformly
  expanding attracting sets. That is, given a continuous
  function we consider its space average with respect to a
  physical measure and compare this with the time averages
  along orbits of the map, showing that the Lebesgue measure
  of the set of points whose time averages stay away from
  the space average tends to zero exponentially fast with
  the number of iterates involved.  As easy by-products we
  deduce escape rates from subsets of the basins of physical
  measures for these types of maps. The rates of decay are
  naturally related to the metric entropy and pressure
  function of the system with respect to a family of
  equilibrium states. \textbf{The corrections added to the
    published version of this text appear in bold; see last
    section for a list of changes}.
\end{abstract}

\maketitle

\tableofcontents

%%%%%%%%%%%%%%%%%%%%%%%%%%%%%%%%%%%%%%%%%%%%%%%%%%%%%%%%%%%%a

\section{Introduction}
\label{sec:intro}

Smooth Ergodic Theory provides asymptotic information on the
behavior of a dynamical system, given by a smooth
transformation, when times goes to infinity. This
statistical approach to Dynamics has provided valuable
insights into many phenomena: from the remarkable result of
Jakobson \cite{Ja81} (see also \cite{BC85,BC91}) showing the
existence of many (positive Lebesgue measure of) parameters
$a\in (0,2)$ for which the corresponding map of the
quadratic family $x\mapsto a-x^2$ has positive Lyapunov
exponent along almost every orbit; to the study of higher
dimensional systems: related ideas provided the first clue
to the nature of the H\'enon attractor \cite{BC91,MV93} or
the existence of robust classes of maps which are not
uniformly expanding but exhibit several distinct positive
Lyapunov exponents \cite{Vi97}, and enabled one to
understand the statistical properties of these and other
classes of systems
\cite{PS82,BeY92,Yo98,BeY99,Al00,BoV00,ABV00}.

% a different set of ideas
% in higher dimensions provided the first clue to the nature
% of the H\'enon attractor \cite{BC91,MV93} or the existence
% of robust classes of maps which are not uniformly expanding
% but exhibit several distinct positive Lyapunov exponents
% \cite{Vi97}, to the study of the statistical properties of
% these and other classes of systems
% \cite{PS82,BeY92,Yo98,BeY99,Al00,BoV00,ABV00}.

The basic ideas can be traced back to the Boltzmann Ergodic
Hypothesis from Statistical Mechanics which was the main
motivation behind the celebrated Birkhoff's Ergodic Theorem
ensuring the equality between temporal and spatial averages
with respect to a (ergodic) probability measure $\mu$ invariant under
a measurable transformation $f:M\to M$ of a compact manifold
$M$, i.e. for every continuous map $\vfi:M\to\RR$ we have
\begin{align}
  \label{e-birkhoff}
  \lim_{n\to+\infty}\frac1n\sum_{j=0}^{n-1}
  \vfi\big( f^j(x) \big) = \int \vfi \, d\mu
\end{align}
for $\mu$ almost every point $x\in M$. Defining $B(\mu)$,
the \emph{ergodic basin of $\mu$}, to be the set of points
for which \eqref{e-birkhoff} holds for every continuous
function $\vfi$, the Ergodic Theorem says that $\mu\big(
B(\mu)\big)=1$ for all ergodic $f$-invariant probability
measures $\mu$. Since ergodic measures can be, for instance,
Dirac masses concentrated on periodic orbits, the Ergodic
Theorem in itself does not always provide information about
the asymptotic behavior of ``big'' subsets of points. The
notion of ``big'' can arguably be taken as meaning ``having
positive Lebesgue measure (or positive volume)'', since such
sets are in principle ``observable sets'' when interpreting
$f:M\to M$ as a model of physical, biological or economic
phenomena.  Correspondingly invariant probability measures
$\mu$ for which $B(\mu)$ has positive volume are called
\emph{physical} (or Sinai-Ruelle-Bowen) measures.

This kind of measures was first constructed for (uniformly)
hyperbolic diffeomorphisms by Sinai, Ruelle and Bowen
\cite{Si72,Ru76,BR75}. Such measures for non-uniformly
hyperbolic maps where obtained more recently
\cite{PS82,BeY92,BeY93,Al00}.

We say that a local diffeomorphism $f$ of a compact manifold
is (uniformly) \emph{expanding} if there exists $n\ge1$ such
that for all $x$ and every  tangent vector $v$ at $x$
\begin{align}
  \label{eq:unifexp}
  \| Df^n(x)v\|\ge 2\|v\|.
\end{align}
For diffeomorphisms of compact manifolds, the notion of
\emph{hyperbolicity} requires the existence of two
complementary directions given by two (continuous)
subbundles $E$ and $F$ of the tangent bundle admitting
$n\ge1$ such that for all points $x$ and  tangent
vectors $(u,v)\in E_x\oplus F_x$
\begin{align}
\label{eq:unifhyp}
  \| Df^n(x)u\|\le \frac12 \|u\|
\quad\mbox{and}\quad
\| Df^n(x)v\|\ge 2\|v\|.
\end{align}
The probabilistic properties of physical measures are an
object of intense study, see e.g.
\cite{BR75,Yo98,BeY99,AA03,alves-araujo2004,alves-luzzatto-pinheiro2005,gouezel,arbieto-matheus2006}.
The leitmotif is that the sequence $\{ \vfi\circ
f^n\}_{n\ge0}$ should behave like an i.i.d. random variable,
at least asymptotically.

Here we are concerned with the rate of convergence of the
time averages \eqref{e-birkhoff}
% in the basin of physical measures 
for non-uniformly expanding maps (NUE) and partially
hyperbolic non-uniformly expanding diffeomorphisms (PHNUE),
where condition~\eqref{eq:unifexp} and the right hand side
condition of~\eqref{eq:unifhyp} are replaced by the
following asymptotic ones
\begin{description}
\item[NUE] for Lebesgue almost every point $x$ there exists
  $n=n(x)\ge1$ such that $\| Df^n(x)v\|\ge 2\|v\|$ for all
   vector $v\in T_x M$;
\item[PHNUE] for Lebesgue almost all points $x$ there exists
  $n=n(x)\ge1$ such that $\| Df^n(x)v\|\ge 2\|v\|$ for all
   vector $v\in F_x$.
\end{description}
We note that if conditions NUE or PHNUE
hold for \emph{every point} then the system is uniformly
expanding or uniformly hyperbolic
\cite{sturman-stark2000,alves-araujo-saussol}.  We also
consider transformations which are diffeomorphisms outside a
``small'' (zero volume) set of singular or critical points
such that the orbits of Lebesgue almost all points have slow
recurrence near this singular set. For more details see the
statement of results below.

The question of the speed of convergence to equilibrium
arises naturally from so-called thermodynamical formalism of
(uniformly) hyperbolic diffeomorphisms, borrowed from
statistical mechanics by Ruelle, Sinai and Bowen (among
others, see
e.g. \cite{Bo75,Ru89,ruelle2004,ellis06,BDV2004})
through the dictionary between one-dimensional lattices and
(uniformly) expanding maps (Gibbs distributions and
equilibrium states in particular) provided by the existence
of a finite Markov partition for the latter systems. Indeed
chaotic dynamics is associated with loss of memory and
creation of information (two views of the same phenomenon)
as the system evolves. These notions are formalized in a
variety of ways, from \emph{entropy}, the exponential rate
of creation of information; to \emph{decay of correlations},
which measures the speed the system ``forgets'' its initial
state; through \emph{large deviations} results, which
measure how fast the system approaches a state of
equilibrium after evolving from almost every initial
state. However, even with abundance of positive Lyapunov
exponents, which is the essential content of the non-uniform
expansion/hyperbolicity conditions above, extending this
theory from uniform to the non-uniform hyperbolic setting
demands considering (if one is optimistic), through the
dictionary already mentioned, Markov partitions with
infinitely many symbols leading to a thermodynamical
formalism of gases with infinitely many states, a hard
subject not yet well understood (see
e.g. \cite{buzzi-sarig2003,arbieto-matheus2006} for recent
developments).

Assuming conditions NUE or PHNUE we are able to extend some
of the large deviation results for uniformly hyperbolic
system in \cite{Yo90,kifer90} (see also
\cite{De89,DenKass2001} for sharp estimates though a
different approach) and strengthen, in a definite sense, the
idea that non-uniformly hyperbolic systems are
\emph{chaotic}: they satisfy a version of the classical
large deviation results for i.i.d.  random variables.  More
precisely, if we set $\delta>0$ as an acceptable error
margin and consider
\[
B_n=\Big\{
x\in M : \Big|
\frac1n\sum_{j=0}^{n-1}
  \vfi\big( f^j(x) \big) 
-
\int\vfi\,d\mu
\Big|>\delta
\Big\}
\]
then we are able to ascertain whether the Lebesgue measure
of $B_n$ decays to zero exponentially fast, i.e. weather
there are constants $C,\xi>0$ such that
\begin{align}
\label{e-decaymeasure}
\Leb\big( B_n \big) \le C e^{-\xi n}
\quad\mbox{for all}\quad n\ge1.
\end{align}
The values of $C,\xi>0$ above depend on $\delta,\vfi$ and on
global invariants for the map $f$ such as the metric entropy
and the pressure function of $f$ with respect to some
equilibrium measures, as detailed in the next section.

%Assuming conditions \textbf{NUE} or \textbf{PHNUE} 
We are able to obtain large deviation rates as
in~\eqref{e-decaymeasure} for non-uniformly expanding local
diffeomorphisms and also for endomorphisms and maps with
non-flat singularities and critical points under a condition
on the rate of approximation of most orbits to the
critical/singular set. In particular we are able to obtain
an exponential decay rate as above for piecewise expanding
maps with infinitely many smoothness domains, for quadratic
maps corresponding to a positive Lebesgue measure subset of
parameters and for a class of maps with infinitely many
critical points. Moreover we also obtain the same kind of
rates for partially hyperbolic attracting sets with a
non-uniformly expanding direction.

%%%%%%%%%%%%%%%%%%%%%%%%%%%%%%%%%%%%%%%%%%%%%%%%%%%%%%%%%%%%a

\subsection{Statement of the results}
\label{sec:statement-results}

We denote by $\|\cdot\|$ a Riemannian norm on the compact
boundaryless manifold $M$, by $d$ the induced distance and
by $\Leb$ a Riemannian volume form, which we call
\emph{Lebesgue measure} or \emph{volume} and assume to be
normalized: $\Leb(M)=1$.

We start by describing one of the class of maps that we are going
to consider. Let $f: M\to M$ be a map of the compact
manifold $M$ which is a $C^2$ local diffeomorphism outside a
set $\cS\subset M$ with zero Lebesgue measure. We assume
that $f$ {\em behaves like a power of the distance} close to
$\cS$: there are constants $B>1$ and $\beta>0$ for which

 \begin{itemize}
 \item[(S1)]
\hspace{.1cm}$\displaystyle{\frac{1}{B} d(x,\cS)^{\beta}\leq
\frac{\|Df(x)v\|}{\|v\|}\leq B d(x,\cS)^{-\beta}}$;
 \item[(S2)]
\hspace{.1cm}$\displaystyle{\left|\log\|Df(x)^{-1}\|-
\log\|Df(y)^{-1}\|\:\right|\leq
B\frac{ d(x,y)}{d(x,\cS)^{\beta}}}$;
 \item[(S3)]
\hspace{.1cm}$\displaystyle{\left|\log|\det Df(x)^{-1}|-
\log|\det Df(y)^{-1}|\:\right|\leq
B\frac{d(x,y)}{d(x,\cS)^{\beta}}}$;
 \end{itemize}
 for every $x,y\in M\setminus \cS$ with $d(x,y)<d(x,\cS)/2$
 and $v\in T_x M\setminus\{0\}$.  The singular set $\cS$ may
 be thought of as containing those points $x$ where $Df(x)$
 is either not defined or else is non-invertible. Note in
 particular that $\cS$ contains the set $\cC$ of critical
 points of $f$, i.e. the set of points (which may be empty)
 where $Df(x)$ is not invertible. We refer to this kind of
 singular sets as \emph{non-flat} since conditions (S1) to
 (S3) above are natural generalizations to arbitrary
 dimensions of the notion of non-flat critical point from
 one-dimensional dynamics, see e.g.\cite{MS93}.

In what follows we write $S_n\vfi(x)$ for
$\sum_{i=0}^{n-1}\vfi(f^i(x))$ and a function
$\vfi:M\to\RR$. 
We say that $f$ as above is \emph{non-uniformly expanding}
if there exists $c>0$ such that
\begin{align}
  \label{eq:NUE}
  \limsup_{n\to+\infty}\frac1n S_n\psi(x) \le -c
\quad\mbox{where}\quad
\psi(x)=\log\big\| Df(x)^{-1} \big\|,
\end{align}
for Lebesgue almost every $x\in M$.
We need to control the rate of approximation of most orbits
to the singular set. We say that $f$ has \emph{slow
  recurrence to the singular set $\cS$} if for every
$\epsilon>0$ there exists $\delta>0$ such that
\begin{align}
  \label{eq:SlowApprox}
  \limsup_{n\to\infty} \frac1n S_n\Delta_\delta(x) < \epsilon
\quad\mbox{with}\quad
\Delta_\delta(x)=\big| \log d_{\delta}(x,\cS) \big|
\end{align}
for Lebesgue almost every $x\in M$, where for any given
$\delta>0$ we define the \emph{smooth  $\delta$-truncated distance}
from $x\in M$ to $\cS$ by
\[
d_\delta(x,\cS)= \xi_\delta\big( d(x,\cS) \big)\cdot
d(x,\cS) + 1- \xi_\delta\big( d(x,\cS) \big)
\]
where $\xi_\delta:\RR\to[0,1]$ is a standard $C^\infty$
auxiliary function satisfying
\[
\xi_\delta(x)=1 \mbox{  if  } |x|\le\delta \mbox{  and  }
\xi_\delta(x)=0 \mbox{  if  } |x|\ge2\delta.
\]
Observe that $\Delta_\delta$ is non-negative and continuous
away from $\cS$ and identically zero $2\delta$-away from
$\cS$.

These notions where presented in~\cite{ABV00} for higher
dimensional maps abstracted from similar notions from
one-dimensional maps \cite{MS93} and previous work on maps
with singularities \cite{KS86}, and in \cite{ABV00,Ze03} the
following result on existence of finitely many physical
measures was obtained.

\begin{theorem}
  \label{thm:abv}
  Let $f:M\to M$ be a $C^2$ local diffeomorphism outside a
  non-flat singular set $\cS$. Assume that $f$ is
  non-uniformly expanding with slow recurrence to $\cS$.
  Then there are finitely many physical (or
  \emph{Sinai-Ruelle-Bowen}) measures $\mu_1,\dots,\mu_k$
  whose basins cover the manifold Lebesgue almost
  everywhere, that is $ B(\mu_1)\cup\dots\cup B(\mu_k) =
  M,\quad \Leb-\bmod0.  $
\end{theorem}

We say that $f$ is a \emph{regular map} if $f_*\Leb\ll\Leb$,
that is, if $E\subset M$ is such that $\Leb(E)=0$, then
$\Leb\big(f^{-1}(E)\big)=0$.  We denote by $\M_f$ the family
of all invariant probability measures with respect to $f$,
by $\M_f^e$ the family of all \emph{ergodic} $f$-invariant
probability measures, and define
\[
B(x,n,\epsilon)=\left\{
y\in M : d\big(f^i(x),f^i(y)\big) < \epsilon, i=0,\dots, n-1
\right\}
\]
the $(n,\epsilon)$-dynamical ball around $x\in M$.  Large
deviation statements are usually related to \emph{local
  entropies} which originated from the works of Shannon,
McMillan and Breiman
\cite{shannon1948,mcmillan1953,breiman1957} and can be
succinctly expressed as follows on a metric space after the
work of Brin and Katok \cite{brin-katok1983}. For any
finite Borel measure $m$ on $M$ define its local entropy at
$x$ to be
\[
h_m(f)(x) = \lim_{\epsilon\to0} \limsup_{n\to\infty}
-\frac1n\log m\Big( B(x,n,\epsilon) \Big).
\]
In \cite{brin-katok1983} it is proved that this limit exists
$m$-almost everywhere whenever $m$ is a $f$-invariant
probability measure. The metric (or measure-theoretic)
entropy of the map $f$ is then defined to be the
non-negative number
\begin{align*}
  h_m(f) = \int h_m(f)(x) \, dm (x).
\end{align*}
Moreover the function $h_m(f)(x)$ is $f$-invariant, so it is
almost everywhere constant if $m$ is $f$-ergodic.

We will be interested in the case $m=$Lebesgue measure
(volume) on $M$, which is usually \emph{not} an invariant
measure in our setting and for $\nu\in\M_f$ we consider
\[
h_m(f,\nu) = \nu-\esup h_m(f).
\]
% Let $\V_f$ be the family of  continuous
% functions $\xi:M\to\RR$ admitting a constant $C>0$ such that
% for Lebesgue almost every $x$ and there exists an integer
% sequence $(n_k)_{k\ge1}$ satisfying for all $k\ge1$
% \[
% \Leb\Big( B(x,n_k,\de_1) \Big) \le C e^{-S_{n_k}\xi(x)}.
% \]
Note that given $\nu\in\M_f$ the value of $h_\nu(f)$ is not
at all related to $h_{\Leb}(f,\nu)$, unless both measures
coincide and $\nu\in\M_f^e$, in which case
$h_\nu(f,\nu)=h_\nu(f)$.

\begin{maintheorem}
\label{mthm:largedeviation}
Let $f:M\to M$ be a regular $C^{1+\alpha}$ local
diffeomorphism outside a non-flat singular set $\cS$, for
some $\alpha\in(0,1)$.  Assume that $f$ is non-uniformly
expanding with slow recurrence to $\cS$. Then writing
$J=\log|\det Df|$, given $c\in\RR$ and a continuous
function $\vfi:M\to\RR$
\begin{enumerate}
\item if $h_\mathrm{top}(f)<\infty$, then
\begin{align*}
\liminf_{n\to+\infty}\frac1n
\log &\Leb\Big(
\big\{
x\in M : \frac1n S_n\vfi(x) > c
\big\}
\Big)
\\
&\ge
\sup
\left\{
h_\nu(f) - h_{\Leb}(f,\nu) :
\nu\in\M_f^e, \int\vfi \,d\nu > c 
\right\};
\end{align*}

\item if $\cS=\emptyset$ ($f$ is a local diffeomorphism)
  then
\begin{align*}
\limsup_{n\to+\infty}&\frac1n
\log \Leb\Big(
\big\{
x\in M : \frac1n S_n\vfi(x) \ge c\big\}
\\
&\le
\sup\left\{
h_\nu(f) - \int J \, d\nu :
\nu\in\M_f, \int\vfi \,d\nu \ge c
\right\}.
\end{align*}

\item in general for any given $\eta>0$ there exists
  $\epsilon,\delta>0$ such that
\begin{align*}
\limsup_{n\to+\infty}&\frac1n
\log \Leb\Big(
\big\{
x\in M : \frac1n S_n\vfi(x) \ge c
\mbox{  and  }
\frac1n S_n\Delta_\delta(x) \le \epsilon
\big\}
\Big)
\\
&\le
\eta+\sup
\left\{
h_\nu(f) - \int J \, d\nu :
\nu\in\M_f, \int\vfi \,d\nu \ge c
\mbox{  and  } \Delta_\delta\in L^1(\nu) 
\right\}.
\end{align*}
\end{enumerate}
\end{maintheorem}

We say that a measure $\nu\in\M_f$ is an
\emph{equilibrium state for $f$ with respect to $J$} (or
just an \emph{equilibrium state} in what follows) if
\[
h_\nu(f) = \nu(J) = \int J \, d\nu.
\]
As the above statement shows, equilibrium states are
involved in the determination of the asymptotic rates of
deviation.  Given $\epsilon,\delta>0$ we write
$\EE=\EE_{\epsilon,\delta}$ for the family of all
equilibrium states $\mu$ of $f$ with respect to $J$ such
that $\mu(\Delta_\delta)\le\epsilon$ and, given a continuous
$\vfi:M\to\RR$, we define $\EE(\vfi)=\{ \nu(\vfi) :
\nu\in\EE\}$.

\begin{remark}
  \label{rmk:different-sups}
  Note that the expressions obtained in items (1) and (2) of
  the statement of Theorem~\ref{mthm:largedeviation} are \emph{not
  comparable} since the supremum is taken over all invariant
  measures in item (2), while we consider only ergodic
  invariant measures in item (1).
\end{remark}

\textbf{We say that $\mu\in\M_f$ is a \emph{weak expanding
    measure} if the subset of points
 weakly satisfying~\eqref{eq:NUE} has full $\mu$-measure, that is}
\begin{align*}
  \mu\left\{x\in M: \limsup_{n\to+\infty}
  \frac1nS_n\psi(x)\le0\right\}=1.
\end{align*}

We are able to deduce that the supremum in the statement of
Theorem~\ref{mthm:largedeviation} is strictly negative for
non-uniformly expanding maps with slow recurrence to the
singular set \textbf{such that all equilibrium states $\EE$
  are expanding measures}.

\begin{maintheorem}
  \label{mthm:supnegative}
  Let $f:M\to M$ be a local diffeomorphism outside a
  non-flat singular set $\cS$ which is non-uniformly
  expanding, has slow recurrence to $\cS$ \textbf{and every
    element in $\EE$ is weak expanding}.  For $\omega>0$ and a
  continuous function $\vfi:M\to\RR$ % admitting $\mu\in\EE$
%   such that $\mu(\vfi)$ is an isolated point in $\EE(\vfi)$
%   and for
% \[
% 0<\omega<\inf\{
%   |\mu(\vfi)-\eta(\vfi)|: \eta\in\EE\setminus\{\mu\} \},
%   \]
there exists $\epsilon,\delta>0$ arbitrarily close to $0$
such that, writing
\[
A_n=\{x\in M: \frac1nS_n\Delta_\delta(x)\le\epsilon\}
\] 
and
\begin{align}
\label{e-defBn}
B_n=\left\{
x\in M : 
\inf\big\{\big|
\frac1n S_n\vfi(x) - \eta(\vfi)
\big| : \eta\in\EE \big\}
> \omega
\right\}
\end{align}
we get
\begin{align}
\label{e-deviationA}
\limsup_{n\to+\infty}\frac1n
\log \Leb\big(A_n\cap B_n\big) <0.
\end{align}
\end{maintheorem}

Clearly if $\cS=\emptyset$ ($f$ is a local diffeomorphism)
then $A_n=M$ and we obtain an asymptotic large deviation
rate for the sets $B_n$. Otherwise to get a similar upper
bound for $\Leb(B_n)$ we need an extra assumption on the
decay of the measure of the \emph{tail sets} $M\setminus
A_n$.

\begin{maincorollary}
  \label{mcor:supnegative1}
  In the setting of Theorem~\ref{mthm:supnegative} with
  $\cS\neq\emptyset$, if $f$ also satisfies
  \begin{align}
    \label{eq:exptail}
\limsup_{n\to\infty}\frac1n\log\Leb(M\setminus A_n ) <0    
  \end{align}
then we have also
\[
\limsup_{n\to\infty}\frac1n\log\Leb( B_n ) <0.
\]
\end{maincorollary}

\begin{remark}
  \label{rmk:integrability} Observe that if $\mu$ is a
  $f$-ergodic absolutely continuous probability measure
  whose support is the entire manifold, then the slow
  recurrence condition \eqref{eq:SlowApprox} is the same as
  saying that $\log d(x,\cS)$ is $\mu$-integrable.

  Note that for any $C^2$ \emph{endomorphism} $f$ (i.e.  the
  singular set $\cS$ of $f$ coincides with the critical set
  $\cC$ of $f$) we have $|\log d(x,\cC)|\ge\Delta_\delta(x)$
  and, as shown in \cite{Li98}, the function $\log d(x,\cC)$
  is $\mu$-integrable for every $f$-invariant probability
  measure. However we need to deal with families of
  invariant probability measures for which $\log d(x,\cC)$
  is \emph{uniformly integrable} so that the proofs of
  Theorems~\ref{mthm:largedeviation}
  and~\ref{mthm:supnegative} can be carried out with our
  arguments. This is why we need the % asymptotic condition on
  % the measures of the
  sets $A_n$ in the previous statements.
  To the best of our knowledge no such general integrability
  result for $\log d(x,\cS)$ exists with respect to
  invariant probability measures for maps with non-flat
  singularities.
\end{remark}

%%%%%%%%%%%%%%%%%%%%%%%%%%%%%%%%%%%%%%%%%%%%%%%%%%%%%%%%%%%%a

\subsection{Partially hyperbolic diffeomorphisms}
\label{sec:part-hyperb-diff}

Let now $f:M\to M$ be a $C^2$ diffeomorphism. We say that a
compact $f$-invariant set $\Lambda$ is an \emph{attracting
  set} if it admits a \emph{trapping region}, that is, an
open neighborhood $U\supset\Lambda$ such that
$\overline{f(U)}\subset U$ and $\Lambda=\cap_{n\ge0}
f^n(U)$. Note that we may have $\Lambda=U=M$ (where $M$
is connected).

As shown in \cite{Yo90}, for every attracting set $\Lambda$
and every physical probability measure $\nu$ supported in
$\Lambda$, given $\delta>0$ and a continuous $\vfi:\overline
U\to\RR$ we have
\begin{align*}
  \liminf_{n\to\infty}
  &\frac1n\log\Leb\left\{
      \Big|
      \frac1n S_n\vfi - \int\vfi\,d\mu
      \Big| > \delta
    \right\}
    \ge
    \\
    &\sup\left\{
      h_\nu(f)-\int\Sigma^+\,d\nu:
      \nu\in\M_f^e, 
      \Big|
      \int\vfi\,d\nu - \int\vfi\, d\mu
      \Big|\ge\delta
    \right\}.
\end{align*}
Here $\Sigma^+$ denotes the sum of the positive Lyapunov
exponents at a given point of $M$. Recall that Ruelle's
Inequality $h_\mu(f)\le\int \Sigma^+ \,d\mu$ is true of
every $C^1$-diffeomorphism \cite{Ru78}.

An attracting set $\Lambda$ is \emph{partially hyperbolic}
(see e.g. \cite{PS82,BDV2004}) if there
exists a continuous splitting $E\oplus F$ of the tangent
bundle of $M$ over $\Lambda$ along two complementary vector
subbundles satisfying
\begin{itemize}
\item $Df$-invariance: $Df(E_x)=E_{f(x)}$ and
  $Df(F_x)=F_{f(x)}$ for all $x\in\Lambda$;
\item domination: there exists $n\ge1$ such that
  \begin{align*}
    \|Df^n\mid E_x\|\cdot\|(Df^n\mid F_x)^{-1}\| 
    \le\frac12
    \quad\mbox{for all}\quad x\in\Lambda;
  \end{align*}
\item $E$ is uniformly contracting: there is $n\ge1$
  such that $\|Df^n\mid E_x\|\le\frac12$ for all $x\in\Lambda$.
\end{itemize}

In this setting we denote by $J$ the logarithm of the
Jacobian along the centre-unstable direction
$J(x)=\log\big|\det Df\mid F_x\big|$ and by $\EE$ the family
of all \emph{equilibrium states} with respect to $J$,
i.e. the set of all $f$-invariant probability measures $\nu$
such that $h_\nu(f)=\nu(J)$.

We will assume further that the $F$ direction only has
positive Lyapunov exponents in the following sense,
introduced in \cite{ABV00}. We say that a partially
hyperbolic attractor with trapping region $U$ is
\emph{non-uniformly expanding} if there exists $c>0$ such
that
\begin{align*}
\limsup_{n\to\infty}\frac1n\sum_{j=0}^{n-1}
\log\big\| (Df\mid F_{f^j(x)})^{-1}\big\| \le -c
\end{align*}
for Lebesgue almost every point $x\in U$.  In
\cite{ABV00} the following was obtained.

\begin{theorem}
  \label{thm:abvdiffeo}
  Let $\Lambda$ be a partially hyperbolic non-uniformly
  expanding attracting set for a $C^2$ diffeomorphism $f$
  with trapping region $U$. Then there are finitely many
  equilibrium states which are physical measures supported
  in $\Lambda$, and whose basins cover $U$ except for a
  subset of zero Lebesgue measure.
\end{theorem}

\textbf{We say that a measure $\mu\in\M_f$ supported in $U$
  is \emph{weak expanding} if the subset of points
  satisfying a weak non-uniformly expanding condition has
  full $\mu$-measure, that is}
\begin{align*}
  \mu\left\{x\in U: \limsup_{n\to\infty}\frac1n\sum_{j=0}^{n-1}
\log\big\| (Df\mid F_{f^j(x)})^{-1}\big\|\le0 \right\}=1.
\end{align*}
We are able to obtain an upper bound entirely analogous to
item 2 of Theorem~\ref{mthm:largedeviation} replacing $M$ by
the points in the trapping region $U$ of a partially
hyperbolic non-uniformly expanding attracting set $\Lambda$
for a $C^2$ diffeomorphism.  Then for the same kind of
attracting sets we obtain an upper bound for the subset
corresponding to \eqref{e-defBn}.

\begin{maintheorem}
  \label{mthm:phdiff<0}
  Let $f:M\to M$ be a $C^2$ diffeomorphism exhibiting a
  partially hyperbolic non-uniformly expanding attracting
  set $\Lambda$ with isolating neighborhood
  $U\supset\Lambda$ \textbf{such that every measure in $\EE$
    is weak expanding}.  Given $\omega>0$ and a continuous
  $\vfi:\overline U\to\RR$, define
\begin{align*}
  B_n=\left\{
x\in U : 
\inf\big\{\big|
\frac1n S_n\vfi(x) - \eta(\vfi)
\big| : \eta\in\EE \big\}
> \omega
\right\}.
\end{align*}
Then
\begin{align*}
  % \label{e-lddiffeo}
    \limsup_{n\to\infty} \frac1n\log\Leb(B_n)<0.
\end{align*}
\end{maintheorem}

%%%%%%%%%%%%%%%%%%%%%%%%%%%%%%%%%%%%%%%%%%%%%%%%%%%%%%%%%%%%a

\subsection{Escape rates}
\label{sec:escape-rates}

Using the estimates obtained above and the observation that
for any compact subset $K$ and a given $\epsilon>0$ we can
find an open set $W\supset K$ and a continuous function
$\vfi:M\to\RR$ such that
\begin{itemize}
\item $\Leb(W\setminus K) < \epsilon$;
\item $0\le\vfi\le1$, $\vfi\mid K\equiv 1$ and $\vfi\mid
  (M\setminus W) \equiv 0$,
\end{itemize}
we see that for $n\ge1$
\begin{align}
  \label{eq:escape}
\left\{
x\in K : f(x)\in K, \dots, f^{n-1}(x)\in K 
\right\}
\subset
\left\{
x\in M : \frac1n S_n\vfi(x)\ge1 
\right\}
\end{align}
and so we get the following (recall the definition of $A_n$
in the statement of Theorem~\ref{mthm:supnegative}).

\begin{maincorollary}
  \label{mcor:escaperate}
  Let $f:M\to M$ be a local diffeomorphism outside a
  non-flat singular set $\cS$ which is non-uniformly
  expanding, has slow recurrence to $\cS$ \textbf{and every
    measure in $\EE$ is weak expanding}.  Let $K$ be a compact
  subset such that $\mu(K)<1$ for all $\mu$ in the
  weak$^*$-closure $\overline\EE$ of $\EE$.  Then for a pair
  $\epsilon,\delta>0$ close to $0$
\[
\limsup_{n\to+\infty}\frac1n
\log \Leb\Big(
\left\{
x\in K\cap A_n : f^j(x)\in K, j=1,\dots,n-1
\right\}
\Big) < 0.
\]
Moreover if $\limsup_{n\to\infty}
\frac1n\log\Leb(M\setminus A_n)<0$ then
\[
\limsup_{n\to+\infty}\frac1n
\log \Leb\Big(
\left\{
x\in K, f(x)\in K, \dots,f^{n-1}(x)\in K
\right\}
\Big)<0.
\]
\end{maincorollary}

In the setting of a partially hyperbolic non-uniformly
expanding attracting set we get, using the same
reasoning as above

\begin{maincorollary}
  \label{mcor:escaperatediffeo}
  Let $f:M\to M$ be a diffeomorphism and $\Lambda$ a
  partially hyperbolic non-uniformly expanding attracting
  set with isolating neighborhood $U$ \textbf{such that
    every measure in $\EE$ is weak expanding}.  Let $K\subset U$
  be a compact subset such that $\mu(K)<1$ for all $\mu$ in
  the weak$^*$-closure $\overline\EE$ of $\EE$.  Then
\[
\limsup_{n\to+\infty}\frac1n
\log \Leb\Big(
\left\{
x\in K, f(x)\in K, \dots,f^{n-1}(x)\in K
\right\}
\Big)<0.
\]
\end{maincorollary}

%%%%%%%%%%%%%%%%%%%%%%%%%%%%%%%%%%%%%%%%%%%%%%%%%%%%%%%%%%%%a

\subsection{Comments and organization of the paper}
\label{sec:organization-paper}

All the arguments use in fact that $f$ is $C^1$ and that its
derivative $Df$ is $\alpha$-H\"older continuous with respect
to the fixed Riemannian norm on $M$, so that all we need is
$f$ to be a $C^{1+\alpha}$ local diffeomorphism outside the
singular set, for some $\alpha\in(0,1)$.

The difficulties we face when considering transformations
which are not uniformly hyperbolic and present singularities
are related to the construction of the measures $\nu$,
appearing in the supremum at item (1) of the statement of
Theorem~\ref{mthm:largedeviation}, as a weak$^*$ limit of
discrete measures which converge to an invariant measure and
are supported on the set one wishes to control. Since we
need to take weak$^*$ limits of measures against
discontinuous test functions, the main body of work in this
paper is to provide sufficient estimates for convergence
imposing some conditions on the dynamics of the maps involved.

The existence of a lower bound for the large deviation rate
with the same expression as in items (2) and (3) of the
statement of Theorem \ref{mthm:largedeviation} depends on
the existence and uniqueness of equilibrium states (the
reader should see \cite{kifer90} for precise statements and
also for counter-examples when uniqueness is not
satisfied). However existence and uniqueness of equilibrium
states for non-uniformly expanding maps is still an open
problem for most potentials in spite of recent progress in
this direction by several authors, see
e.g. \cite{OliVi2005,AMSV06,arbieto-matheus2006}.

Recently Pinheiro \cite{Pinheiro05} has extended the
statement of Theorem~\ref{thm:abv} replacing the limsup in
condition \eqref{eq:NUE} by liminf, keeping the same
conclusions involving the existence of finitely many
physical measures and of a positive density of hyperbolic
times Lebesgue almost everywhere. Hence our statements are
automatically valid in this more general setting.

In what follows, we start by presenting some non-trivial
classes of maps to which our results are applicable, in
Section~\ref{sec:examples-application}. In
Section~\ref{sec:hyperbolic-times} we present preliminary
technical results to be used in the following sections.
Theorem~\ref{mthm:largedeviation} is then proved in
Subsection~\ref{sec:upper-bound-large} for local
diffeomorphisms, in Subsection~\ref{sec:upper-bound-diffeo}
for partially hyperbolic non-uniformly expanding
diffeomorphisms and in Subsection~\ref{sec:upper-bound-with}
for maps with singularities or criticalities.  In
Section~\ref{sec:striclly-negat-upper} we deduce
Theorem~\ref{mthm:supnegative} from
Theorem~\ref{mthm:largedeviation}, first for local
diffeomorphisms and for the partially hyperbolic case in
Subsection~\ref{sec:local-diff-case}, and then with
singularities or criticalities in
Subsection~\ref{sec:case-with-sing}, together with an
extension of Ruelle's Inequality to maps with non-flat
singularities in Subsection~\ref{sec:entr-form-piec}.

\subsubsection*{Acknowledgements}

We are thankful to M. Viana (IMPA) for valuable comments and
suggestions during the elaboration of this
text. \textbf{Later, P. Varandas (UFBA) pointed out to the
  authors several issues with the proofs that prompted the
  preparation of this corrected version.} The authors are
also indebted to the fine scientific environment and access
to the superb mathematical library of IMPA during the
preparation of the earliest versions of the manuscript.

%%%%%%%%%%%%%%%%%%%%%%%%%%%%%%%%%%%%%%%%%%%%%%%%%%%%%%%%%%%%a

\section{Examples of application}
\label{sec:examples-application}

Here we show that there are many examples of maps in the
conditions of Theorem~\ref{mthm:supnegative},
Corollary~\ref{mcor:supnegative1} or
Theorem~\ref{mthm:phdiff<0}.

\subsection{Quadratic maps and infinite-modal maps}
\label{sec:quadratic-maps}

In~\cite{ArPa04} the following $C^\infty$ family of maps of
$I=[-1,1]$ with infinitely many critical points was considered:
\begin{align*}
%\label{e2.4,8}
f_{\mu}(z)=\left\{
\begin{array}{ll}
f(z)+\mu & \mbox{   for   } z\in (0,\epsilon]\\
f(z)-\mu & \mbox{   for   } z\in [-\epsilon,0)
\end{array}
\right.
\end{align*}
where $f:I\to I$ is an expanding extension of
\begin{align*}
%\label{e2.3}
\hat{f}:[-\epsilon,\epsilon]\to[-1,1],
\quad
\hat{f}(z)=\left\{
\begin{array}{ll}
az^{\alpha}\sin(\beta\log(1/z))) & \mbox{   if   }z>0\\
-a|z|^{\alpha}\sin(\beta\log(1/|z|))) & \mbox{   if   }z<0,
\end{array}
\right.
\end{align*}
to $I$ (i.e. $|f'|\gg1$ on
$I\setminus[-\epsilon,\epsilon]$), with $a>0$,
$0<\alpha<1, \, \beta>0$ and $\epsilon>0$. It was shown that
there exists a positive Lebesgue measure subset $P$ of
parameters in $(-\epsilon,\epsilon)$ such that for
$\mu\in P$ the map $f_\mu$ is non-uniformly expanding and
has slow recurrence to the non-flat infinite and denumerable
singular set. Moreover for the same parameters the decay
rate of the tail set is exponential, i.e. \eqref{eq:exptail}
is true.  \textbf{If all equilibrium states with respect to
  $-\log|f'|$ are weak expanding}, then $f_\mu$ for
$\mu\in P$ is in the setting of Corollaries
\ref{mcor:supnegative1} and \ref{mcor:escaperate}.

Analogous results hold for the quadratic family
$Q_a(x)=a-x^2$ (and also for general $C^2$ unimodal
families), so that Corollaries \ref{mcor:supnegative1} and
\ref{mcor:escaperate} apply to quadratic maps for a positive
Lebesgue measure subset of parameters \textbf{since all
  invariant measures are weak expanding in this setting;
  see~\cite{Przyty93}}. % Observe that since
% $Q_a$ is a endomorphism of a compact non-empty interval, it
% is in the setting of Subsection~\ref{sec:isolated}.

%%%%%%%%%%%%%%%%%%%%%%%%%%%%%%%%%%%%%%%%%%%%%%%%%%%%%%

\subsection{Piecewise smooth one-dimensional expanding maps}
\label{sec:piecewise-smooth-one}

Let $f:I\to I$ be a map admitting a sequence
$\cS=\{a_n, n\ge1\}\subset I=[-1,1]$ such that for every
connected component $G$ of $I\setminus\cS$ we have that
$f\mid G$ is $C^1$ diffeomorphism with its image and
\textbf{there exists $n\in\ZZ^+$ so that $|Df^n(x)|>1$ for
  all $x\in I\setminus\cup_{i=0}^{n-1}f^{-i}\cS$}. Assume
that $\cS$ is a non-flat singular set for $f$ and that $f$
admits a absolutely continuous ergodic invariant probability
measure $\mu$ with positive Lyapunov exponent and such that
$\log d(x,\cS)$ is $\mu$-integrable and $\supp\mu=I$. Then
$f$ is in the setting of Theorem~\ref{mthm:supnegative}
\textbf{since all invariant measures are weak expanding in
  this case}.

Examples of this kind of maps are the Gauss map
\cite{Vi97b}, and transitive piecewise one dimensional maps
satisfying the conditions in \cite{Ry83} (see also
\cite{Vi97b}), that is there exists $\kappa>0$ such that for
every connected component $G$ of $I\setminus\cS$ we also have
\begin{align*}
  \var_G\frac1{|f'|}\le \kappa \cdot\sup_G\frac1{|f'|}
\quad\mbox{and}\quad
\sum_{G} \sup_G\frac1{|f'|} \le \kappa.
\end{align*}
More concrete examples are Lorenz-like maps
\cite{LY73,Vi97b} (even with criticalities \cite{LV00}) and
the maps introduced by Rovella \cite{Ro93,mtz001}.

A proof of the exponential decay of the tail set for this
class of maps is not available in the literature to the best
of our knowledge but can be done as an application of the
technique of exclusion of parameters introduced in
\cite{BC85} (the details will appear in forthcoming work
\cite{araujo2006a}), so that Corollaries
\ref{mcor:supnegative1} and \ref{mcor:escaperate} also hold
for this type of maps.

% \begin{remark}
%   \label{rmk:existsingular}
%   From \cite{Le81} it is known for piecewise $C^2$ maps $f$
%   of the interval or the circle with finitely many domains
%   of smoothness that, if an invariant probability measure
%   $\mu$ has positive entropy, then $\mu\in\EE$ if, and only
%   if, $\mu$ is absolutely continuous with respect to
%   Lebesgue measure. This enables us to argue as in
%   Subsection~\ref{sec:isolated} and conclude that, in this
%   setting, if $f$ is non-uniformly expanding and has slow
%   recurrence to the singular set, then there exists isolated
%   elements in $\EE$.
% \end{remark}

%%%%%%%%%%%%%%%%%%%%%%%%%%%%%%%%%%%%%%%%%%%%%%%%%%%%%%

\subsection{Non-uniformly expanding local diffeomorphisms}
\label{sec:non-unif-expand}

Consider a local diffeomorphism $f:M\to M$, so that
$\cS=\emptyset$, which satisfies
\begin{itemize}
\item $\|(Df)^{-1}\|\le1$ and
\item $K_1=\{ x\in M : \|Df(x)^{-1}\|=1\}$ is finite.
\end{itemize}
Then by the results in~\cite{ArTah} we have that such $f$
has a finite set $\EE$ of ergodic equilibrium states for
$\phi$ \textbf{all of which are weak expanding measures}.
Hence in this case Theorem~\ref{mthm:supnegative} holds for
every continuous function $\vfi:M\to\RR$.

\subsection{Viana maps}
\label{sec:viana-maps}

The following family of endomorphisms of the cylinder
was introduced by Viana in \cite{Vi97}.  Let $a_0\in(1,2)$
be such that the critical point $x=0$ is preperiodic for
the quadratic map $Q(x)=a_0-x^2$. Let $\sS^1=\RR/\ZZ$ and
$b:\sS^1\rightarrow \RR$ be a Morse function, for instance
$b(s)=\sin(2\pi s)$. For fixed small $\alpha>0$, consider
\begin{align*}
\begin{array}{rccc} \hat f: & \sS^1\times\RR
&\longrightarrow & \sS^1\times \RR\\
 & (s, x) &\longmapsto & \big(\hat g(s),\hat q(s,x)\big)
\end{array}
\end{align*}
where $\hat g$ is the uniformly expanding map of the circle
defined by $\hat{g}(s)=d\cdot s$ (mod $\ZZ$) for some $d\ge16$,
and $\hat q(s,x)=a(s)-x^2$ with $a(s)=a_0+\alpha b(s)$.  For
$\alpha>0$ small enough there exists an interval $I\subset (-2,2)$
such that $\hat f(S^1\times I)$ is contained in the interior
of $S^1\times I$. Hence any map $f$ sufficiently $C^0$ close to
$\hat f$ has $S^1\times I$ as a
forward invariant region. We consider from here on these
maps $f$ close to $\hat f$ restricted to $\sS^1\times I$.

In \cite{Vi97,Al00,AA03} a $C^3$ neighborhood $\U$ of $\hat
f$ was studied and it was proved that every $f\in\U$ is
non-uniformly expanding and has slow recurrence to the
non-flat critical set $\cC$. The arguments in \cite{Vi97}
where extended in \cite{buzzi-sester-tsujii} to encompass
the weaker condition $d\ge2$ on the expansion of $\hat g$,
providing the same properties for a $C^\infty$-neighborhood
$\tilde\U$ of $\hat f$.

\sout{Hence, each $f\in\U$ or $f\in\tilde\U$ is in the setting of
Theorem~\ref{mthm:supnegative}.} Results in
\cite{alves-luzzatto-pinheiro2005,gouezel} show that the
tail set decays at least sub-exponentially fast, which is
not enough to ensure that
Corollaries~\ref{mcor:supnegative1} and
\ref{mcor:escaperate} are true for the maps in $\U\cup
\tilde\U$. It is conjectured that the tail set indeed decays
exponentially fast and with a uniform rate for all maps in
$\U\cup\tilde\U$.

%%%%%%%%%%%%%%%%%%%%%%%%%%%%%%%%%%%%%%%%%%%%%%%%%%%%%%%%%%

\subsection{Partially hyperbolic non-uniformly expanding
  diffeomorphisms}
\label{sec:open-class-part}

We sketch the construction of a robust class of partially
hyperbolic non-uniformly expanding diffeomorphisms, taking
$U$ equal to $M$, following~\cite{ABV00}. This construction
is closely related to the $C^1$ open classes of transitive
non-Anosov diffeomorphisms presented in
\cite[Section~6]{BoV00}, as well as other robust examples from
\cite{Man87}.

Start with a linear Anosov diffeomorphism $\hat f$ on the
$d$-dimension\-al torus $M=\TT^d$, $d\ge 3$. Write
$TM=E\oplus F$ the corresponding hyperbolic decomposition of
the tangent bundle with $\dim F\ge2$. Let $V$ be a small
closed domain in $M$ for which there exist unit open cubes
$K^0$ and $K^1$ in $\RR^d$ such that $V \subset \pi(K^0)$
and $\hat f(V)\subset \pi(K^1)$, where $\pi:\RR^d\to \TT^d$
is the canonical projection. Let now $f$ be a diffeomorphism
on $\TT^d$ such that
 \begin{enumerate}
 \item[(A)] $f$ admits invariant cone fields $C_E$ and
   $C_F$, with small width $a>0$ and containing,
   respectively, the stable bundle $E$ and the unstable
   bundle $F$ of $\hat f$;
 \item[(B)] $f$ is \emph{partially hyperbolic and volume
     expanding along the center-unstable direction}:
 there is $\sigma_1>1$ so that
 \begin{align*}
|\det(Df\mid T_x\D_F)| > \sigma_1 
\quad\mbox{and}\quad 
\|Df\mid T_x\D_E\| < \sigma_1^{-1}   
 \end{align*}
 for any $x\in M$ and any disks $\D_F$, $\D_E$ tangent to
 $C_F$, $C_E$, respectively (see
 Subsection~\ref{sec:cover-part-hyperb} for more on
 invariant cone fields and disks tangent to cone fields in
 this setting).
 \item[(C)] $f$ is $C^1$-close to $\hat f$ in
the complement of $V$, so that there exists $\sigma_2<1$
satisfying
 $$
 \|(Df \mid T_x \D_F)^{-1}\| < \sigma_2
\quad\mbox{and}\quad \|Df \mid T_x \D_E\| < \sigma_2
 $$
for any $x\in (M\setminus V)$ and any disks $\D_F$,
$\D_E$ tangent to $C_F$, $C_E$,
respectively. Moreover $f(V)$ is also contained in
the projection of a unit open cube.
\item[(D)] there exist some small $\delta_0>0$ satisfying $$
  \|(Df \mid T_x\D_F)^{-1}\| < 1+\delta_0 $$
  for any
  $x\in V$ and any disk $\D_F$ tangent to $C_F$.
\end{enumerate}

If $\tilde f$ is a torus diffeomorphism satisfying (A), (B),
(D), and coinciding with $\hat f$ outside $V$, then any map
$f$ in a $C^1$ neighborhood of $\tilde f$ satisfies all the
previous conditions. Results in~\cite[Appendix]{ABV00} show
in particular that for any $f$ satisfying (A)--(D) there
exist $c_u>0$ such that $f$ is partially hyperbolic and
non-uniformly expanding along its center-unstable direction,
as defined in
Subsection~\ref{sec:part-hyperb-diff}.  Hence on a
small $C^2$ neighborhood $\U$ of $\tilde f$ every
diffeomorphism $f\in\U$ satisfies all the conditions of
Theorem~\ref{mthm:phdiff<0} \textbf{if  all
  equilibrium states with respect to the central-unstable
  Jacobian have only non-negative Lyapunov exponents along the
  central-unstable direction. This can be achieved by
  certain $C^1$ perturbations of a linear Anosov
  diffeomorphism $\hat f$.}

%%%%%%%%%%%%%%%%%%%%%%%%%%%%%%%%%%%%%%%%%%%%%%%%%%%%%%%%%%

\section{Hyperbolic times}
\label{sec:hyperbolic-times}

The main technical tool used in the study of non-uniformly
expanding maps is the notion of hyperbolic times, introduced
in \cite{Pl72,Al00}. We say that $n$ is a
$(\sigma,\delta,b)$-hyperbolic time of $f$ for a point $x$
if the following two conditions hold with $0<\sigma<1$ and
$b,\delta>0$
\begin{align}
    \label{eq:tempo-hip}
\prod_{j=n-k}^{n-1}\big\| Df\big(f^j(x)\big)^{-1}\big\| \le \sigma^k
\quad\mbox{and}\quad
d_\delta \big( f^k(x),\cS \big) \ge e^{-bk}
\end{align}
for all $k=0,\dots,n-1$.

We now outline the properties of these special times. For detailed
proofs see~\cite[Proposition 2.8]{ABV00} and~\cite[Proposition
2.6, Corollary 2.7, Proposition 5.2]{AA03}.

\begin{proposition}
  \label{pr:prophyptimes}
  There are constants $C_1,\delta_1>0$ depending on
  $(\sigma,\de,b)$ and $f$ only such that, if $n$ is
  $(\sigma,\de,b)$-hyperbolic time of $f$ for $x$, then
  there are \emph{hyperbolic preballs} $V_k(x)$ which are
  neighborhoods of $f^{n-k}(x)$, $k=1,\dots, n$, such that
\begin{enumerate} 
\item $f^k\mid V_k(x)$ maps $V_k(x)$ diffeomorphically to
  the ball of radius $\delta_1$ around $f^n(x)$;
\item   for every $1\leq k\leq n$ and $y,z\in V_k(x)$
\[
d\big(f^{n-k}(y),f^{n-k}(z)\big)\le
  \sigma^{k/2}\cdot d\big(f^{n}(y),f^{n}(z)\big);
\]
\item for $y,z\in V_k(x)$
\[
\frac1{C_1}\le
\frac{\big|\det Df^{n-k}(y)\big|}{\big|\det Df^{n-k}(z)\big|} \le
C_1.
\]
\end{enumerate}
\end{proposition}

The following ensures existence of infinitely many
hyperbolic times for Lebesgue almost every point for
non-uniformly expanding maps with slow recurrence to the
singular set. A complete proof can be found in~\cite[Section
5]{ABV00}.

\begin{theorem}
\label{thm:tempos-hip-existem}
Let $f:M\to M$ be a $C^{1+\alpha}$ local diffeomorphism away
from a non-flat singular set $\cS$, for some
$\alpha\in(0,1)$, non-uniformly expanding and with slow
recurrence to $\cS$. Then there are $\sigma\in(0,1)$,
$\delta,b>0$ and there exists
$\theta=\theta(\sigma,\delta,b)>0$ such that $\Leb$-a.e.
$x\in M$ has infinitely many $(\sigma,\de,b)$-hyperbolic
times.  Moreover if we write $0<n_1<n_2<n_2<\dots$ for the
hyperbolic times of $x$ then their asymptotic frequency
satisfies
\[
\liminf_{N\to\infty}\frac{\#\{ k\ge1 : n_k\le
  N\}}{N}\ge\theta
\quad\mbox{for}\quad \Leb\mbox{-a.e.  } x\in M.
\]
\end{theorem}

%%%%%%%%%%%%%%%%%%%%%%%%%%%%%%%%%%%%%%%%%%%%%%%%%%%%%%%%%%

\subsection{Coverings by hyperbolic preballs}
\label{sec:cover-hyperb-pre}

Here we show how to cover a given measurable subset with
hyperbolic preballs, which will enable us to approximate its
Lebesgue measure through the measure of families of
hyperbolic preballs. In turn, the measure of a hyperbolic
preball is related to the Jacobian of the transformations
due to bounded distortion.

\begin{lemma}
  \label{le:meashiptimes}
Let $B\subset M$, $\theta>0$ and $g:M\to M$ be a local
diffeomorphisms outside a non-flat exceptional set
$\cS$ such that $g$ has density $>2\theta$ of hyperbolic
times for every $x\in B$. Then, given any probability
measure $\nu$ on $B$ and any $m\ge1$, there exists $n>m$
such that
\[
\nu\big( \{
x\in B : n \mbox{  is a hyperbolic time of $g$ for  } x
 \} \big) > \frac\theta2.
\]
\end{lemma}

This is \cite[Lemma 4.4]{OliVi2005} easily adapted to our
setting. For completion we include its very short
proof. This lemma shows that we can translate the density of
hyperbolic times into the Lebesgue measure of the set of
points which have a specific (large) hyperbolic time.

\begin{proof}
  Let $H$ be the set of pairs $(x,n)\in B\times\NN$ for
  which $n$ is a hyperbolic time of $g$ for $x$. For each
  $k\ge1$, let $\#_k$ be the normalized counting measure on
  $\{m+1,m+2,\dots,m+k\}$. Our assumption implies that for
  any given $x\in B$ we have for big enough $k\ge1$
\[
h_k(x)=\#_k \big( \pi(H\cap (\{x\}\times\NN) ) \big) > 2\theta,
\]
where $\pi:B\times\NN\to\NN$ is the projection on the second
coordinate.  Given any probability measure $\nu$ on $B$ we
have by Fatou's Lemma
\begin{align*}
  \liminf_{k\to\infty}\int h_k\,d\nu \ge \int
  \liminf_{k\to\infty} h_k \, d\nu \ge 2\theta
\end{align*}
so we may fix $k\ge1$ large enough so that
$\nu(h_k)>\theta$ and find a subset
for $C\subset B$ with $\nu(C)>1/2$ and $h_k(x)\ge \theta/2$
for all $x\in C$. By Fubini's Theorem this means that
\[
(\nu\times\#_k)(H) > \theta\mbox{ and thus } \nu\big(\hat\pi( H\cap(
B\times \{n\} )) \big) > \frac\theta2
\]
for some $m<n\le m+k$, where $\hat\pi:B\times\NN\to B$ is
the projection on the first coordinate. This proves the
lemma.
\end{proof}

\subsubsection{Construction of an adequate initial
  partition}
\label{sec:constr-an-adequate}

Let $f$ be a regular map in the setting of the Main Theorem
with positive density of $(\sigma,\delta)$-hyperbolic times
for Lebesgue almost everywhere \textbf{and
  $\rho:M\setminus\cS\to(0,+\infty)$ a continuous positive
  function possibly unbounded.  Let
  $U_m=\rho^{-1}[e^{-(m+1)},e^{-m})$ and $N\in\ZZ^+$ be such
  that $U_m\neq\emptyset$ and $8e^{-m}<\delta_1$ for $m> N$;
  and also $U_{N}=M\setminus U_{N+1}$. These sets have
  non-empty interior and are relatively compact.}

\textbf{Fix $0<\delta_0<\delta_1$ and let $\B^N$ be a
  finite open cover of $U_{N}$ by $\delta_0$-balls and
  $\B^m$ a finite open cover of $U_m$ by $r_m$-balls, where
  $r_m=\min\{e^{-(n+1)},\delta_0\}/8$ for $n> N$. Since $M$
  is a finite dimensional manifold, we can find such open
  cover with a number $\ell_m$ of $r_m$-balls such that
  $\ell_m\le Cr_m^{-\dim M}$ for all sufficiently small
  $r_m$.}

  \textbf{Let also
  $\B=\cup_{m\ge N}\B^m$ be a countable open cover of
  $M\setminus\cS$ and let us enumerate the elements of
  $\B^N$ first, then $\B^{N+1},\B^{N+2},\dots$ in this
  order, obtaining $\B=\{B_k: k\ge1\}$.
%=\{B(x_i,\delta_1/8), i=1,\dots,l\}
From this we define a \textbf{countable} partition $\cP$ of $M$
such that $\diam \cP(x)<\rho(x)$ following the proof of
\cite[Lemma 13.3]{Man87}}.

\textbf{We start by setting $P_1=B_1\cap U_N$ as the first
  element of the partition $\cP$. Then, assuming that
  $P_1,\dots, P_k$ are already defined we set
  $P_{k+1}=B_{k+1} \setminus(P_1\cup\dots\cup P_k)$ for
  $k+1\le\#\B^N$.}  Note that if $P_k\neq\emptyset$ then
$P_k$ has non-empty interior, diameter smaller than
$\delta_0/4$ and the boundary $\partial P_k$ is a (finite)
union of pieces of boundaries of balls in a Riemannian
manifold, thus has zero Lebesgue measure.  \textbf{This
  provides a partition of $U_N$ whose nonempty atoms we
  include in $\cP$.}

\textbf{We now repeat this procedure for each $m>N$ obtaining a
finite partition of $U_m$ whose nonempty atoms we include in
$\cP$. Note that if $P\in\cP$ and $\emptyset\neq P\subset
U_m$, then $P$ has nonempty interior; $\diam
P\le\min\{e^{-(m+1)},\delta_0\}/8\le\rho(x), \forall x\in P$
and again $\partial P$ is a finite union of pieces of
boundaries of balls in $M$.}

Note that since $f$ is regular
the boundary of $g(P)$ still has zero Lebesgue measure for
every atom $P\in\cP$ and every inverse branch $g$ of $f^n$,
for any $n\ge1$. 

\textbf{Let us choose one interior point in each atom $P\in\cP$
contained in $U_m$ and form the set $\cC_m$ of
representatives of the atoms of $\cP$ in $U_m$; and  let
$d_m=\min\{ d(w,\partial\cP), w\in\cC_m\}>0$ where $m\ge N$ and
$\partial\cP=\cup_{P\in\cP}\partial P$ is the boundary of
$\cP$.}

\begin{proposition}
  \label{pr:zeroboundary}
  Let $(\mu_n)_{n\ge1}$ be a family of Borel probability
  measures on $M$; $\mu$ some weak$^*$ accumulation point of
  the sequence $(\mu_n)$ \textbf{and
    $\rho:M\setminus\cS\to(0,+\infty)$ be a continuous
    $\mu$-integrable function.}  Then, given
  $0<\xi\le\tau$, there exists a partition
  $\cP_{\xi,\tau}$ with \textbf{the same number of
    atoms of $\cP$ in each $U_m, m\ge N$, each atom has
    non-empty interior and zero Lebesgue measure boundary;
    and also}
  \begin{enumerate}
  \item $\mu(\partial\cP_{\xi,\tau})=0$ and
    $\mu_n(\partial\cP_{\xi,\tau})=0$ for all $n\ge1$;
  \item each $P\in\cP_{\xi,\tau}$ contains one, and only one,
    \textbf{element of $\cC=\sum_{m\ge N}\cC_m$}\footnote{We
      write $A+B$ the union of the disjoint subsets $A$ and
      $B$.};
      \item $2\diam\cP_{\xi,\tau}(x)\le\min\{\rho(x),\tau_1\}$ for
    $\Leb$-, $\mu$- and $\mu_n$-a.e. $x$;
  \item for each $P\in\cP_{\xi,\tau}$ there is $Q\in\cP$
    satisfying
    $\Leb(P\triangle Q)<\xi<\tau\cdot\Leb(Q)$;
  \item $H_\mu(\cP_{\xi,\tau})<\infty$.
  \end{enumerate}
\end{proposition}

\begin{proof}
  \textbf{Let $0<\xi<\tau>0$ be given. For each fixed
    $m\ge N$, let us take $0<\gamma_m<\min\{\xi,d_m,r_m^3\}$
    such that for any given $r_m$-ball $B=B(x,r_m)\in\B^m$}
  \begin{align}
\label{eq:smallvolume}
    \Leb\left( B\big(x,r_m + \gamma_m\big) \setminus
B\big(x, r_m\big) \right) < \epsilon_m/\#\B^m
  \end{align}
  where
  $0<\epsilon_m<\min\{\xi,
  \tau\cdot\min\{\Leb(B):B\in\cP, B\subset
  U_m\}\}$; and also for all $n\ge1$
\begin{align}
\label{eq:zeroboundary}
\mu\big( \partial B(x, r_m + \gamma_m) \big)=0
=\mu_n\big( \partial B(x, r_m + \gamma_m)
\big)
\end{align}
and in addition for $a_m=(1-\gamma_m)e^{-m}$
\begin{align}
  \label{eq:levelboundary}
  \Leb(\rho^{-1}(a_m))=0=\mu(\rho^{-1}(a_m))=\mu_n(\rho^{-1}(a_m)).
\end{align}
Such value of $\gamma_m$ exists since the set of such values
so that some of the expressions in~\eqref{eq:zeroboundary}
or~\eqref{eq:levelboundary} is positive for some $B\in\B^m$,
some $m\ge N$ and some $n\ge1$ is denumerable. Thus we may
take $\gamma_m>0$ satisfying~\eqref{eq:zeroboundary}
and~\eqref{eq:levelboundary} arbitrarily close to zero, and
so inequality~\eqref{eq:smallvolume} can also be obtained.

\textbf{We consider now the open cover $\wt{\B}$ of
  $M\setminus\cS$ obtained by replacing $U_m$ by
  $\wt{U_m}=\rho^{-1}[a_{m+1},a_m)$ and each $r_m$-ball of
  $\B^m$ by a concentric $(r_m+\gamma_m)$-ball in
  $\wt{\B^m}$ for each $m\ge N$, and construct the partition
  $\wt{\cP}$ obtained from $\wt{\B}=\sum_{m\ge N}\wt{\B^m}$ by
  the same procedure as before with the same order. Since
  $\gamma_m<\epsilon<d_m$ we obtain
  $d\big(w,\partial \cP_\epsilon)\ge d_m-\gamma_m>0$ for all
  $w\in\cC_m, m\ge N$}.

  This shows that each $w\in\cC$ is contained in some atom
  $P_w$ of $\wt{\cP}$.  Moreover there cannot be distinct
  $w_1,w_2\in\cC$ such that $w_2\in P_{w_1}$, because this
  would mean that for some $m\ge N$ and $B=B(x,r_m)\in\B^m$
  we have $w_2\in B(x,r_m)$, $w_1\not\in B(x,r_m)$ and
  $w_1,w_2\in B(x,r_m+\gamma_m)$, which contradicts the
  choice of $\gamma_m<d_m$.

  Hence, on the one hand, $\#\cP\le\#\wt{\cP}$.  On the
  other hand, let us consider $\{ P_w, w\in\cC\}$. There
  might be other (finitely many) atoms $P$ in $\wt{\cP}$
  and, if so, we join them to some adjacent atom $P_w$
  (meaning $\overline P\cap\overline P_w\neq\emptyset$)
  obtaining a new atom $P\cup P_w$. In this way we obtain a
  partition, which we still denote by $\wt{\cP}$ with as
  many atoms as the elements of $\cC$ and satisfying items
  (1) and (2) of the statement of the lemma.

  \textbf{Finally, for any $w\in\cC$ the corresponding atoms
    $P_w\in\wt{\cP}$ and $Q_w\in\cP$ satisfy
    $P_w\in\wt{\B^m}, Q_w\in\B^m$ for some $m\ge N$ and}
\[
  \Leb\big(P_w\triangle Q_w\big) \le \sum_{B(x,r_m)\in\B^m}
  \Leb\left( B\big(x,r_m + \gamma_m\big)\setminus
    B(x,r_m)\right) < \#\B^m\cdot\epsilon_m/\#\B^m =
  \epsilon_m\le\epsilon
\]
\textbf{and $\diam(P_w)\le 4r_m <\min\{\rho(x),\delta_1\}/2$ for all
  $x\in P_w$. This provides item (3) of the statement of the
  lemma. By the choice of $\epsilon_m$ we also get}
\[
  \Leb\big(P_w\triangle Q_w\big)
  <
  \epsilon_m
  \le
  \tau
  \cdot\min\{\Leb(B):B\in\cP,
  B\subset U_m\}
\le\tau\cdot\Leb(Q_w).
\]
\textbf{This is item (4) of the statement of the lemma. To
  prove that $\wt{\cP}$ has finite entropy, we use
 observe that the number of atoms of $\wt{\cP}$ on each
$U_m$ is bounded by $\ell_m$, and so by construction we obtain}
\begin{align}
  H_\mu(\wt{\cP})
  &=
  \sum_{m\ge N}\sum_{P\in\wt{\cP}: P\subset
  U_m}\hspace{-0.5cm}-\mu(P)\log\mu(P)
  =
  \sum_{m\ge N}\sum_{P\in\wt{\cP}: P\subset
  U_m}\hspace{-0.5cm}\mu(U_n)
  \left(-\frac{\mu(P)}{\mu(U_n)}\log\frac{\mu(P)}{\mu(U_n)}
    -\log\mu(U_n)\right) \nonumber
  \\
  &\le
    \sum_{m\ge N}\mu(U_n)\big(\log\ell_m-\log\mu(U_n)\big)
  \\
  &\le\label{eq:summability}
    \sum_{m\ge N}\left(
    -\mu(U_n)\log\mu(U_n)+\mu(U_n)\log C+\dim(M)(m+1)\mu(U_m)
    \right).
\end{align}
\textbf{Then we deduce the summability of the above series
  using that $\log\rho\in L^1(\mu)$, as follows. On the one
  hand, we note that}
\begin{align}\label{eq:logrho}
  \sum_{m\ge N}(m+1)\mu(U_m)
  =
  1+\sum_{m\ge N}-\log e^{-m}\mu(U_m)
  \le
  1+\int|\log\rho|\,d\mu<\infty.
\end{align}
\textbf{On the other hand, we have}
\begin{lemma}[Lemma 13.2 from~\cite{Man87}]\label{le:summability}
  If $x_n\in(0,1), n\ge1$ and $\sum_{n\ge1}n x_n<\infty$,
  then $\sum_{n\ge1}x_n\log(1/x_n)<\infty$.
\end{lemma}
\textbf{Now setting $x_n=\mu(U_n)$ we have the assumption of
Lemma~\ref{le:summability} from~\eqref{eq:logrho} and so we deduce}
\begin{align*}
  \sum_{n\ge N}-\mu(U_n)\log\mu(U_n)=\sum_{n\ge N}x_n\log(1/x_n)<\infty.
\end{align*}
This completes the proof of the summability
of~\eqref{eq:summability} and with it the proof of
\textbf{Proposition}~\ref{pr:zeroboundary} after setting
$\cP_{\xi,\tau}=\wt{\cP}$.
\end{proof}

\subsubsection{The flexible covering lemma}
\label{sec:flexible-covering-le}

Having this we can now obtain the following flexible
covering lemma.

\begin{lemma}
  \label{le:coverhiptimes}
  Let a measurable set $E\subset M$, $m\ge1$ and
  $\epsilon>0$ be given with $\Leb(E)>0$. Let $\theta>0$ be
  a lower bound for the density of hyperbolic times for
  Lebesgue almost every point. Then there are integers
  $m<n_1<\dots<n_k$ for $k=k(\epsilon)\ge1$ and families
  $\cE_i$ of subsets of $M$, $i=1,\dots, k$ such that
\begin{enumerate}
\item $\cE_1\cup\dots\cup\cE_k$ is a finite
  pairwise disjoint family of subsets of $M$;
\item $n_i$ is a $(\sigma/2,\delta/2)$-hyperbolic time for
  every point in $P$, for every element $P\in\cE_i$, $i=1,\dots,k$;
\item every $P\in\cE_i$ is the preimage of some element
  $Q\in\cP$ under an inverse branch of $f^{n_i}$, $i=1,\dots,k$;
\item there is an open set $U_1\supset E$ containing the
  elements of $\cE_1\cup\dots\cup\cE_k$ with
  $\Leb(U_1\setminus E)<\epsilon$;
\item $\Leb\Big( E\triangle \bigcup_i \cE_i \Big) \le
  \left(1-\frac{\theta}4 \right)^k < \epsilon.$
\end{enumerate}

\end{lemma}

The proof follows \cite[Lemma 8.2]{OliVi2005} closely. We
write $\mC_m$ the set of pairs $(z,n_i)$ where
$f^{n_i}(z)=w\in\cC$ and $z\in P$ for all $P\in\cE_i$ and
$i=1,\dots, k$ (such $z$ exist by item (3) of
Lemma~\ref{le:coverhiptimes}).

\begin{remark}\label{rmk:kdeepsilon}
  Note that  $k$ depends on $\epsilon$ only and
  not on the set $E$.
\end{remark}

\begin{proof}[Proof of Lemma~\ref{le:coverhiptimes}]
  By the non-uniformly expanding assumption on $f$ we know
  that there exists $\theta>0$ such that Lebesgue almost
  every point has density $>\theta$ of hyperbolic times of
  $f$. 

Let $U_1$ be an open set and $K_1$ a compact set such that
$K_1\subset E\subset U_1$ and $\Leb(U_1\setminus
K_1)<\epsilon$ and $\Leb(K_1)>(1/2) \Leb(U_1)$. Using
Lemma~\ref{le:meashiptimes} with $B=K_1$ and
$\nu=\Leb/\Leb(K_1)$ we can find $n_1>m$ such that $e^{-c
  n_1} < d( K_1 , M\setminus U_1)$ and the subset $L_1$ of
points of $K_1$ for which $n_1$ is a hyperbolic time
satisfies $\Leb(L_1)\ge\frac{\theta}2\Leb(K_1)\ge
\frac{\theta}4\Leb(E)$. 

Given $x\in L_1$ let
$g:B(f^{n_1}(x),\delta_1)\to V_{n_1}(x)$ be the inverse
branch of $f^{n_1}\mid V_{n_1}(x)$, recall that $n_1$ is a
hyperbolic time for $x$ and see
Proposition~\ref{pr:prophyptimes}. By the choice of
$\cP$ from \textbf{Proposition}~\ref{pr:zeroboundary} there exists a unique
$P\in\cP$ such that $f^{n_1}(x)\in P$. Let us consider
$g(P)$ and let $\cE_1$ be the family of all such sets
obtained as $g(P)$ which intersect $L_1$, where $g$ is an
inverse branch of $f^{n_1}$ corresponding to a hyperbolic
time and $P$ is an element of $\cP$.
 
Note that the elements of $\cE_1$ are pairwise disjoint
because $\cP$ is a partition.  Moreover by the properties of
hyperbolic times (Proposition~\ref{pr:prophyptimes}) the
diameter of $P\in\cE_1$ is smaller than $e^{-c n_1}$. Hence
the union $E_1$ of all the elements of $\cE_1$ is contained
in $U_1$ and by construction
\[
\Leb(E_1\cap E) \ge \Leb(L_1) \ge \frac{\theta}4 \Leb(E).
\]
Now consider the open set $U_2=U_1\setminus \overline{E_1}$
and set $K_2\subset E\setminus \overline{E_1}$ a compact set
such that $\Leb(K_2)\ge(1/2)\Leb(E\setminus E_1)$. Observe
that $\Leb(\overline{E_1}\setminus E_1)=0$ since
$\partial\cP$ has zero Lebesgue measure and this property is
preserved under backward iteration by the regularity
assumption on $f$. Reasoning as before, we can find
$n_2>n_1$ such that $e^{-c n_2} < d( K_2, M\setminus U_2)$
and a set $L_2\subset K_2$ such that
$\Leb(L_2)\ge(\frac{\theta}2)\Leb(K_2)$ and $n_2$ is a
hyperbolic time for every $x\in L_2$. Let $\cE_2$ be the
family of elements $g(P)$ which intersect $L_2$, where
$P\in\cP$ and $g$ is an inverse branch of $f^{n_1}$
corresponding to a hyperbolic time.

Again $\cE_2$ is a pairwise disjoint family of sets whose
diameters are smaller than $e^{-c n_2}$. Thus their union
$E_2$ is contained in $U_2$. Hence $\cE_1\cup\cE_2$ is also
a pairwise disjoint family and, in addition
\[
\Leb\big(
E_2\cap( E\setminus E_1) 
\big)
\ge \Leb(L_2) 
\ge
\frac{\theta}2 \Leb(K_2)
\ge
\frac{\theta}4 \Leb(E\setminus E_1).
\]
Repeating this procedure we obtain families $\cE_i,
i=1,\dots,k$ of elements of $\cP_{n_i}$ which are pairwise
disjoint and contained in $U_1$, and
\begin{align}
  \label{eq:tends0}
  \Leb\Big(E_{i+1}\cap\big(E\setminus(E_1\cup\dots\cup
  E_i)\big)\Big)\ge
\frac{\theta}4
\Leb\big(E\setminus(E_1\cup\dots\cup
  E_i)\big)
\end{align}
for all $i=1,\dots,k-1$, for some $k\ge1$, where $E_j=\cup
\cE_j$. Hence 
\[
\Leb\Big( \bigcup_{i=1}^k E_i \setminus E\Big)
\le
\Leb( U_1\setminus E) < \epsilon
\]
and (\ref{eq:tends0}) ensures that
\[
\Leb\Big( E \setminus \bigcup_{i=1}^k E_i \Big)
\le 
\left(1-\frac{\theta}4\right)^k\Leb(E).
\]
Therefore we can find $k\ge1$ such that $\Leb\big(E\triangle
\cup_{i=1}^k \cE_i\big)<\epsilon$, as stated.
\end{proof}

\begin{remark}
  \label{rmk:partitionperturbation}
  Note that the construction proving
  Lemma~\ref{le:coverhiptimes} gives a finite sequence of
  hyperbolic times, open sets $U_1,\dots, U_k$ and closed
  sets $\overline E_1,\dots,\overline E_k$. Having these we
  can find small enough $\delta>\epsilon>0$, replace $\cP$
  in the proof of Lemma~\ref{le:coverhiptimes} by any
  partition $\cP_{\epsilon,\delta}$ obtained as in
  \textbf{Proposition}~\ref{pr:zeroboundary} (by slightly
  modifying $\cP$), and use the same inverse branches to
  obtain families $\cE_i^\prime$ of preballs such that
\[
\Leb\left(\Big( \bigcup_i \cE_i\Big)
\triangle \Big(\bigcup_i
  \cE_i^\prime \Big)
\right) 
\le 
\sum_i C_1\delta \Leb(\cE_i)
<
C_1\delta\Leb\Big( \bigcup_i \cE_i\Big)
\le 
C_1\delta
\]
where $C_1$ is the volume distortion constant (see
Proposition~\ref{pr:prophyptimes}). Hence, after the
modification of the initial partition, we get
\[
\Leb\big(E\triangle\bigcup_i
  \cE_i^\prime\big) < \epsilon+C_1\delta < (1+C_1)\delta
\]
since $\epsilon<\delta$.  Moreover the set $\mC_m$ is
unaffected since $\cC$ is fixed and the inverse branches
are kept.
\end{remark}

%%%%%%%%%%%%%%%%%%%%%%%%%%%%%%%%%%%%%%%%%%%%%%%%%%%%%%%

\subsection{The partially hyperbolic setting}
\label{sec:cover-part-hyperb}

Here we state the main results needed to obtain an extension
of the covering Lemma~\ref{le:coverhiptimes} to the setting
of partially hyperbolic non-uniformly expanding attracting
sets. As we indicate along the way, the proofs of most of
them can be found in~\cite{ABV00}.

\subsubsection{Stable/Unstable cone fields}

Let $\Lambda$ be a partially hyperbolic and non-uniformly
expanding attracting set for a $C^2$ diffeomorphism $f:M\to
M$ with a trapping region $U\subset M$.  The existence of
the dominated splitting $E\oplus F$ of $T_\Lambda M$ ensures
the existence of a continuous extension $\tilde E\oplus
\tilde F$ of $E\oplus F$ to a neighborhood of $\Lambda$,
which we assume without loss to be $U$, and of the following
cone fields:
\begin{description}
\item[stable cones] $\EE^a_x=\{(u,v)\in \tilde E(x)\oplus
\tilde  F(x): \|v\|\le a\cdot\|u\|\}$;
\item[unstable cones] $\FF^b_x=\{(u,v)\in \tilde E(x)\oplus
\tilde  F(x): \|u\|\le b\cdot\|v\|\}$;
\end{description}
for all $x\in U$ and $a,b\in(0,1)$, which are $Df$-invariant
in the following sense (see e.g.
\cite[Appendix C]{BDV2004})
\begin{itemize}
\item if $x,f^{-1}(x)\in U$, then $Df^{-1}(\EE^{ a}_x)\subset
  \EE^{\lambda a}_{f^{-1}(x)}$;
\item if $x,f(x)\in U$, then $Df(\FF^b_x)\subset
  \FF^{\lambda b}_{f(x)}$;
\end{itemize}
for some $\lambda\in(0,1)$.  Continuity enables us to
unambiguously denote $d_E=\dim (\tilde E)$ and $d_F=\dim
(\tilde F)$, so that $d=d_E+d_F=\dim(M)$, and domination
guarantees that the angles between the $\tilde E$ and
$\tilde F$ directions are bounded from below away from zero
at every point.

\subsubsection{Hyperbolic times}

In this setting, given $\sigma>1$ we say that $n$ is a
$\sigma$-hyperbolic time for $x\in U$ if
\begin{align*}
  \prod_{j=n-k+1}^{n}\big\| (Df \mid
  F_{f^{j}(x)})^{-1}\big\| \le \sigma^k \qquad\text{for all
    $1\le k \le n$.}
\end{align*}

\begin{remark}
  \label{rmk:hyptimesdiffeo}
This definition of hyperbolic time is entirely analogous to
the one given in the local diffeomorphisms setting except
that we restrict the derivatives to the $F$-direction. Hence
the statement and proof of Lemma~\ref{le:meashiptimes} carry
over without change.
\end{remark}

\subsubsection{$E$-disks and $F$-disks}

Let us fix the unit balls of dimensions $d_E,d_F$
\[
\BB_E=\{w\in\RR^{d_E}: \|w\|_2\le1\}\qand
\BB_F=\{w\in\RR^{d_F}: \|w\|_2\le1\}
\]
where $\|\cdot\|_2$ is the standard Euclidean norm on the
corresponding Euclidean space.  We say that a $C^{1+\alpha}$
embedding $\De:\BB_E\to M$ (respectively $\De:\BB_F\to M$)
is a \emph{$E$-disk} (resp. \emph{$F$-disk}) if the image of
$D\De(w)$ is contained in $\EE^a_{\De(w)}$ for all
$w\in\BB_E$ (resp. $D\De(w)(\RR^{d_F})\subset
\FF_{\De(w)}^b$ for every $w\in\BB_F$), where
$\alpha\in(0,1)$ if fixed.

\subsubsection{Curvature of $E$- and $F$-disks at hyperbolic times}

Let $r_0>0$ be an injectivity radius of the exponential map on
$M$, that is $\exp_x:B(x,r_0)\to M$ is a
diffeomorphism onto its image
$G(x,r_0)=\exp_x\big(B(x,r_0)\big)$, where $B(x,r_0)=\{v\in
T_X M: \|v\|<r_0\}$ is the $r_0$-neighborhood of $0$ in $T_x
M$.  By the continuity of the splitting $E\oplus F$ and the
cone fields we can choose $0<r<\min\{r_0,\delta_1/4\}$ such
that for every $x\in\Lambda$ the subspace $E_x$ is contained
in all the images of the cone field $\EE_x^a$ under the
exponential map $\exp_x$ and analogously for the
complementary direction, that is for every $y\in
G(x,r)\cap\Lambda$ we have
\begin{align}
  \label{e-cone}
  E_x\subset D (\exp_x^{-1})\big( \EE_y^a\big)
  \quad\mbox{and}\quad
  F_x\subset D (\exp_x^{-1})\big( \FF_y^b\big).
\end{align}
This ensures that every $F$-disk (respectively every
$E$-disk) $\De$ is such that its image on $B(x,r)$ given
by $\exp_x^{-1}\big(\De\cap G(x,r)\big)$ is transversal
to the direction of $E_x$ (resp. $F_x$).

The ``curvature'' of $E$- and $F$-disks can be determined by
the notion of H\"older variation of the tangent bundle as
follows.  
%  We write $V_x=B(x,r_0)$ in what follows.  We
% are going to identify $V_x$ through the local chart
% $\exp_x^{-1}$ with the neighborhood $U_x=\exp_x (V_x)$ of
% the origin in $T_x M$, and we also identify $x$ with the
% origin in $T_x M$. In this way we get that $E(x)$ (resp.
% $F(x)$) is contained in $\EE_y^a$ (resp. $\FF_y^b$) for all
% $y\in U_x$, reducing $\de_0$ if needed, and the intersection
% of $F(x)$ with $\EE_y^a$ (and the intersection of $E(x)$
% with $\FF_y^b$) is the zero vector.

We write $\De$ also for the image of the respective
embedding for every $E$- or $F$-disk.  Hence if $\De$ is a
$E$-disk and $y=\De(w)$ for some $w\in\BB_E$, then the
tangent space of $\De$ at $y$ is the graph of a linear map
$A_x(y):T_x\De \to F(x)$ for $w\in\De^{-1}(V_x)$ (here
$T_x\De=D\De(x)(\RR^{d_E})$). The same happens locally for a
$F$-disk exchanging the roles of the bundles $E$ and $F$
above.

The domination condition on the splitting $E\oplus F$
ensures the existence of $\zeta\in(0,1)$ such that for some
$n\ge1$ and all $x\in\Lambda$
\begin{align*}
  \|Df^n\mid E_x\|\cdot\|(Df^n\mid F_x)^{-1}\|^{1+\zeta} 
    \le\frac34.
\end{align*}
Given $C>0$ we say that the \textit{tangent bundle of $\De$
  is $(C,\zeta)$-H\"older} if
\begin{equation}
  \label{eq:holder-bundle}
  \| A_x(y) \| \le C \dist_\De(x,y)^\zeta \quad\mbox{for all}\quad
y\in G(x,r)\cap \De \qand x\in U,
\end{equation}
where $\dist_\De(x,y)$ is \textit{the distance along $\De$}
defined by the length of the shortest smooth curve from $x$
to $y$ inside $\De$ calculated with respect to the
Riemannian norm $\|\cdot\|$ induced on $TM$.

For a $E$- or $F$-disk $\De\subset U$ we define
\begin{equation}
  \label{eq:curvature-definition}
  \kappa(\De)=\inf\{C>0 : T\De \mbox{ is }(C,\zeta)\mbox{-H\"older} \}.
\end{equation}
The proof of the following result can be found
in~\cite[Subsection 2.1]{ABV00}. The basic ingredients are
the cone invariance and dominated decomposition properties
for $f$.

\begin{proposition}
\label{pr:bounded-curvature}
There is $C_2>0$ such that given a $F$-disk $\De\subset U$
  \begin{enumerate}
  \item there exists $n_1\in\NN$ such that
  $\kappa(f^n(\De))\le C_2$ for all $n\ge
  n_1$;
  \item if $\kappa(\De)\le C_2$ then
  $\kappa(f^n(\De))\le C_2$ for all
  $n\ge0$;
  \item in particular, if $\De$ is as in the previous item,
  then
    $$
    J_n: f^n(\De) \ni x
    \mapsto \log | \det (Df \mid T_x (f^n(\De)) |
    $$
    is $(L_1,\zeta)$-H\"older continuous with $L_1>0$
    depending only on $C_2$ and $f$, for every $n\ge1$.
  \end{enumerate}
\end{proposition}

% \begin{proof}
%   See \cite[Proposition~2.2]{ABV00} and
%   \cite[Corollary~2.4]{ABV00}.
% \end{proof}

\subsubsection{Distortion bounds}

The following uniform backward contraction and distortion
bounds are proved in \cite[Lemma 2.7, Proposition
2.8]{ABV00}.

\begin{proposition}
\label{p.contraction}
There exist $C_3,\delta_1>0$ depending only on
$f,\sigma$ such that, given any $F$-disk $\Delta\subset U$,
$x\in \Delta$, and $n \ge 1$ a $\sigma$-hyperbolic time for
$x$,
\begin{enumerate}
\item $\dist_{f^{n-k}(D)}(f^{n-k}(y),f^{n-k}(x)) \le \sigma^{k/2}
\dist_{f^n(D)}(f^n(y),f^n(x)), $ for all $y\in \Delta$ satisfying
$\dist(f^n(x),f^n(y))\le\delta_1$;
\item if $\kappa(\Delta) \le C_2$ then
  \begin{align*}
    \frac{1}{C_3} 
    \le 
    \frac{|\det Df^{n} \mid T_y \Delta|}
    {|\det Df^{n} \mid T_x \Delta|}
    \le C_3
  \end{align*}
  for every $y\in \Delta$ such that $\dist(f^n(y),f^n(x))\le
  \delta_1$.
\end{enumerate}
\end{proposition}

\subsubsection{The initial partition and the covering lemma}

Now we consider the following rectangle
\begin{align*}
  \hat R(x,s)=\{(u,v)\in T_x M:\|u\|<s,\, \|v\|< s,\, u\in
  E_x,\, v\in F_x\}
\end{align*}
where $s$ is chosen so that $\hat R(x,s)\subset B_x(r)$ for
all $x\in\Lambda$. This defines an open cover
$\{\exp_x\big(\hat R(x,s)\big)\}_{x\in\Lambda}$ of $\Lambda$
which admits a finite subcover denoted by
$\R=\{R_1=R(x_1,s),\dots, R_h=R(x_h,s)\}$. This finite cover will
define the initial partition $\cP$ 
% However we need to ``refine the
% initial partition along the $E$ direction'' when dealing
% with big hyperbolic times to ensure that the diameter of the
% atoms decreases to zero fast enough. The initial partition
% $\cP_1$ 
given by
\begin{align*}
\cP=\{R_1,M\setminus
R_1\}\vee\dots\vee\{R_h,M\setminus R_h\}.  
\end{align*}
% Consider a partition $\D_{n,i}$ of the disk $\{u\in E_{x_i}
% : \|u\|<s\}$ defined through a finite open cover
% $B_{i,1},\dots, B_{i,h_i}$ by $r_n$-balls as
% \begin{align*}
% \D_{n,i}=\{B_{i,1},M\setminus
% B_{i,1}\}\vee\dots\vee\{B_{i,h_i},M\setminus B_{i,h_i}\}  
% \end{align*}
% for each $i=1,\dots,l$.  Then consider the refinement of
% $\cP_1$ induced by the $\D_{n,i}$ as follows
% \begin{align*}
%   \cP_n=\{ \exp_{x_i}\big( D\times \tilde D_i\big):
%   D\in\D_{n,i}, i=1,\dots,l\}
% \end{align*}
% where $\tilde D_i=\{v\in F_{x_i}: \| v \| < r \}$, see
% Figure~\ref{fig:cover-refined}. Observe that
% $\cP_1\le\cP_2\le\cP_3\le\dots$ by construction.
We may assume without loss that $\Leb(\partial\cP)=0$ by
slightly changing the initial cover.  We choose an interior
point in each element of $\cP$ which together define the
set $\cC$.
% \begin{figure}[htbp]
%   \centering
%   \includegraphics[width=10cm,height=4cm]{fig-cover-refined.eps}
%   \caption{\label{fig:cover-refined}The original cover with
%     solid lines and the refined covers in the $E$ direction
%     in dashed lines.}
% \end{figure}

Now we adapt the covering Lemma~\ref{le:coverhiptimes} to
the setting of partially hyperbolic non-uniformly expanding
attracting sets as follows.

\begin{lemma}
  \label{le:coverpartialhyp}
  Let a measurable set $E\subset U$, $m\ge1$ and
  $\epsilon>0$ be given. Let $\theta>0$ be a lower bound for
  the density of hyperbolic times for Lebesgue almost
  every point on $U$. Then there are 
  integers $m<n_1<\dots<n_k$ for $k=k(\epsilon)\ge1$, and
  families $\cE_i$ of subsets of $M$, $i=1,\dots, k$ such
  that
\begin{enumerate}
\item $\cE_1\cup\dots\cup\cE_k$ is a finite
  family of subsets of $M$ and each $\cE_i$ is a pairwise
  disjoint family;
\item $n_i$ is a $(\sigma/2,\delta/2)$-hyperbolic time for
  every point in $P$, for every element $P\in\cE_i$, $i=1,\dots,k$;
\item every $P\in\cE_i$ is the preimage of some element
  $Q\in\cP$ under $f^{-n_i}$, $i=1,\dots,k$;
% \item there is an open set $U_1\supset E$ containing the
%   elements of $\cE_1\cup\dots\cup\cE_k$ with
%   $\Leb(U_1\setminus E)<\epsilon$;
\item $
\Leb\Big( E\setminus \bigcup_i \cE_i \Big) \le
  \left(1-\frac{\theta}4 \right)^k < \epsilon.$
\end{enumerate}
\end{lemma}

\begin{proof}
  Let $E\subset U$, $\epsilon>0$ and $m\ge1$ be given. Set
  $\nu=\Leb/\Leb(E)$ and apply Lemma~\ref{le:meashiptimes}
  with $B=E$ to obtain $n_1>m$ and $L_1\subset E$ such that
  $n_1$ is a hyperbolic time for every point $x\in L_1$ and
  $\Leb(L_1)\ge\frac{\theta}2 \Leb(E)$.

  Given $x\in L_1$ let $P_x$ be the unique element of the
  partition $f^{-n_1}\cP$ which contains $x$ (recall that
  $f$ is a diffeomorphism). Define $\cE_1=\{P_x: x\in
  L_1\}$. Then $\cE_1$ is a finite pairwise disjoint family
  of preimages of elements of $\cP$ corresponding to a
  hyperbolic time $n_1$. If $E_1$ is the union of the
  elements of $\cE_1$, then
  \begin{align*}
    \Leb(E_1\cap E) \ge \Leb(L_1) \ge \frac{\theta}2 \Leb(E).
  \end{align*}
  Now consider $\hat E_2=E\setminus\overline E_1$. If
  $\Leb(\hat E_2)<\epsilon$ then we are done, since then
  $\Leb(E\setminus E_1)<\epsilon$ because
  $\Leb(\partial\cE_1)=0$ as $f$ is regular map. Otherwise
  use again Lemma~\ref{le:meashiptimes} to find $n_2>n_1$
  and $L_2\subset \hat E_2$ such that $n_2$ is a hyperbolic
  time for all points of $L_2$ and
  $\Leb(L_2)\ge\frac{\theta}2 \Leb(\hat E_2)$.

  Let $\cE_2$ be the family of all elements of the partition
  $f^{-n_2}\cP$ which intersect $\hat E_2$. Then $\cE_2$ is
  a pairwise disjoint family and the union $E_2$ of its
  elements satisfies
  \begin{align*}
    \Leb\big( E_2\cap( E\setminus E_1) \big)
    \ge \Leb(L_2) 
    \ge \frac{\theta}2 \Leb(\hat E_2)
    \ge \frac{\theta}4 \Leb(E\setminus E_1).
  \end{align*}
  Repeating this procedure we get families $\cE_i,
  i=1,\dots,k$ of elements of $f^{-n_i}\cP$ with
  $m<n_1<\dots<n_k$ satisfying the inequality
  \eqref{eq:tends0}. These families satisfy items (1)-(3) by
  construction and item (4) follows by~\eqref{eq:tends0} as
  in the proof of Lemma~\ref{le:coverhiptimes}. This
  concludes the proof.
\end{proof}

Observe that we may apply
\textbf{Proposition}~\ref{pr:zeroboundary} to $\cP$
\textbf{with $\rho=const.$} to ensure that, for a given
denumerable family of $f$-invariant probability measures,
there is a partition $\cP_{\xi,\tau}$ arbitrarily
close to $\cP$, with the same number of elements, such that
the measure of the boundary of the elements of
$\cP_{\xi,\tau}$ is zero with respect to all measures
of the family. Moreover, as in the previous subsection, we
write $\mC_m$ the set of pairs $(z,n_i)$ where
$f^{n_i}(z)=w\in\cC$ and $z\in P$ for all $P\in\cE_i$ and
$i=1,\dots, k$. In addition, we can build the new partition
$\cP_{\xi,\tau}$ in such a way that the sets $\mC_n$
are unchanged.

%%%%%%%%%%%%%%%%%%%%%%%%%%%%%%%%%%%%%%%%%%%%%%%%%%%%%%%%%%%%%%

\subsection{The volume of dynamical balls}
\label{sec:meas-dynam-balls}

Here we show that the volume of dynamical balls on
hyperbolic times is well controlled by $S_n J$, either in
the local diffeomorphism case with or without singularities,
or in the partially hyperbolic case.

\subsubsection{The local diffeomorphism case with singularities}
\label{sec:endomorphism-case-or}

Note that by the properties of bounded distortion of volumes
during hyperbolic times (item 3 of Proposition
\ref{pr:prophyptimes}) we can write, if $n$ is a hyperbolic
time of $f$ for $x\in M$
\begin{align*}
\Leb\big( B(f^k(x),n-k,\delta_1) \big)
&=
\int_{B(f^k(x),n-k,\delta_1)} \frac{dz}{\big|\det Df^{n-k}(z)\big|}
\\
&\le
C_1
\frac{\Leb\big( B(f^n(x),\delta_1)\big)}{\big|\det Df^{n-k}(x)\big|},
\end{align*}
then recalling that $J=\log|\det Df|$ we get
\begin{align*}
\Leb\big( B(f^k(x),n-k,\delta_1) \big)
&\le
C_1 e^{-S_{n-k}J(f^k(x))} \Leb\big( B(f^n(x),\delta_1)
\\
&\le
C_1 e^{-S_{n-k}J(f^k(x))}.
\end{align*}
Observe that by Proposition \ref{pr:prophyptimes} if $n$
is a hyperbolic time of $f$ for $x$ we get due to uniform
backward contraction
\[
S_{n-k}J(f^k(x))=
\log \big|\det Df^{n-k}(x)\big| 
\ge (n-k)\cdot\dim(M)\log\sigma/2 >0
\]
which will be used several times in what follows.

%%%%%%%%%%%%%%%%%%%%%%%%%%%%%%%%%%%%%%%%%%%%%%%%%%%%%%%

\subsubsection{The partially hyperbolic case with non-uniform expansion}
\label{sec:part-hyperb-case}

In the partially hyperbolic and non-uniformly expanding
setting we recall the construction of the cover
$\R=\{R_1,\dots,R_j\}$ and the initial partition $\cP$ from
Subsection~\ref{sec:cover-part-hyperb}. Observe that if we
take $\delta_0$ to be the Lebesgue number of the covering
$\R$ (see e.g. \cite{munkres}), then for all
$0<\delta<\delta_0$ we have for all $x\in U$ and $n\ge1$ a
hyperbolic time for $x$
\begin{align*}
  B(x,n,\delta)\subset f^{-n}\cP(x),
\end{align*}
where $f^{-n}\cP(x)$ denotes the element of $f^{-n}\cP$ which
contains $x$.  To find an upper bound for the volume of this
dynamical ball it is enough to estimate the volume of
$f^{-n}\cP(x)$ when $n$ is a hyperbolic time for $x$.

Let $P\in\cP$ be such that $f^{-n}(P)$ has a positive
Lebesgue measure subset $\tilde P$ of points for which $n$
is a hyperbolic time and choose $h$ such that $R_h\supset
P$.  Let $\tilde Q\in\cP$ be such that $Q=\tilde
Q\cap\tilde P$ has positive Lebesgue measure and choose $l$
such that $R_l\supset Q$.  

We consider the projection of $\hat P=\exp_{x_l}^{-1}(\tilde
P)$ on $E_{x_l}$ parallel to $F_{x_l}$. Its diameter will be
bounded by a constant which is a function of $f$ and $s$ only,
since the number of different $R_l$ is finite.  Projecting
$\hat Q$ on the complementary direction $F_{x_l}$ parallel
to $E_{x_l}$ we may use the backward contraction and bounded
area distortion for hyperbolic times along $F$-disks to
estimate the area along $F$-disks and integrate to deduce a
volume estimate.

\begin{figure}[h]
  \centering
  \includegraphics[width=11cm,height=4cm]{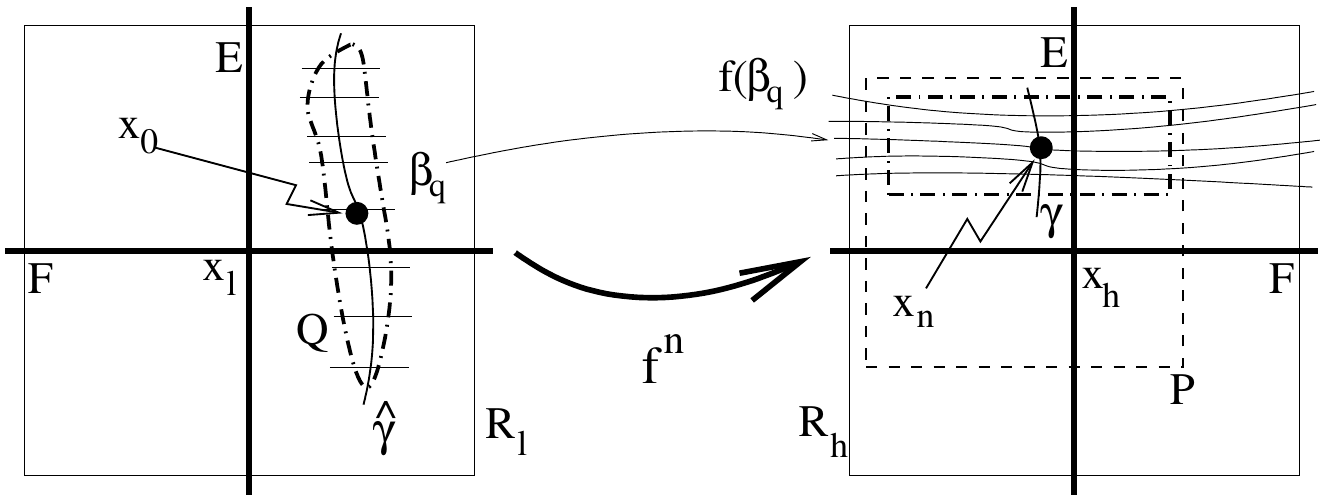}
  \caption{\label{fig-F-disks} The diameter of the elements
    of $\cE_n$ through the use of $E$-disks and images of
    $F$-disks on a hyperbolic time.}
\end{figure}

Indeed, observe that since the $E$ direction is uniformly
contracted by $Df$, if we fix a point $x_0\in Q$, the
corresponding point $x_n=f^n(x_0)\in P\cap f^n(Q)$ and a
$E$-disk $\gamma$ which crosses $R_h$, then the connected
component $\hat\gamma$ of $f^{-n}(\gamma)\cap R_l$
containing $x_0$ is a $E$-disk which also crosses $R_l$.
Moreover distances along $\gamma$ are uniformly expanded by
$f^{-1}$. Thus every point $w_0\in\hat\gamma$ is such that
$w_k=f^k(w_0)$ and $x_k=f^k(x_0)$ satisfy
\begin{align}
\label{eq:unifcontraction}
C\frac{\delta_1}4>  C s \ge \dist(w_0,x_0) 
\ge C \lambda^{-k} \dist(w_k,x_k),
\end{align}
for some constant $C>0$ depending on $f$ only.  Hence if we
take $s$ small enough then we can ensure that $w_k$ is close
enough to $x_k$ for $k=1,\dots, n$ so that $n$ is also an
hyperbolic time for all $w_0\in\hat\gamma$.  Thus we can
consider $F$-disks $\beta_q$ through the points $q$ of $Q$
parallel to $F$, which are transversal to $\hat\gamma$. Then
the images $f^n(\beta_q)$ will be $F$-disks crossing $R_l$
which together cover $P\cap f^n(Q)$, see
Figure~\ref{fig-F-disks}.

The preimages $f^{-n}(P\cap f^n(Q)\cap f^n(\beta_q))$ then
form a cover of $Q$ and these predisks are $F$-disks whose
diameter is smaller than $e^{-c n}$.

% Therefore the diameter of $Q$ is
% smaller than $d'+d''< C e^{-c n}$.
% Note that each element $P\in\cE_i$ is the preimage under
% $f^{-n_i}$ of some element of $\cP_{n_i}$. However since
% $\cP_n\le\cP_m$ for $n\le m$, then each element of
% $\cP_{n_i}$ is a finite pairwise disjoint union of elements
% of $\cP_{n_k}$ for all $i=1,\dots,k$. Thus each $P\in\cE_i$
% is also the union of preimages of elements of $\cP_{n_k}$
% under $f^{-n_i}$. We may without loss redefine $\cE_i$
% replacing each $P\in\cE_i$ by the smaller pieces from
% $f^{-n_i}\cP_{n_k}$, concluding the proof of the new
% statement of item 3.

% The rest of the proof of Lemma~\ref{le:coverhiptimes}
% follows without change in this setting. We define as well
% the sets $\cC_{n_i}$ of interior points of the preballs used
% while constructing the set $\cE_i$, $i=1,\dots,k$.

Using Tonelli's Theorem we can write $\Leb\big( Q \big) =
\int_{\hat\gamma} m\big( Q \cap \beta_q\big) \, dq$ where
$m$ denotes the $d_F$-dimensional Lebesgue measure induced
by $\Leb$ on $F$-disks and $dq$ is Lebesgue measure along
the disk $\hat\gamma$. By the Change of Variables Formula
together with the bounded area distortion along hyperbolic
times in the partially hyperbolic setting given by
Proposition~\ref{p.contraction} we get for each
$q\in\hat\gamma$
\begin{align*}
  m\big( Q \cap \beta_q\big)&=\int_{\beta_q} \chi_Q \,dm
  =\int_{f^{-n}(f^n(\beta_q))}  \chi_Q \,dm
  \\
  &=\int_{f^n(\beta_q)} (\chi_Q\circ f^{-n})\cdot \big|\det
  Df^{-n}\mid f^n(\beta_q)\big| \,dm
  \\
  &=\int_{f^n(\beta_q)} e^{-S_nJ(f^{-n}(z))}
  \chi_{f^n(Q)}(z) \,dm(z)
  \\
  &\le C_3 \cdot e^{-S_nJ(f^{-n}(q))}\cdot m\big(f^n(Q)\cap
  f^n(\beta_q)\big),
\end{align*}
thus $ \Leb\big( Q \big) \le \int_{\hat\gamma} C_3 e^{-S_n
  J(q)} m\big(f^n(Q)\cap f^n(\beta_q)\big) \, dq$.
But by~\eqref{eq:unifcontraction} we see that every
$q\in\hat\gamma\cap Q$ satisfies 
\begin{align*}
  d(f^k(q),f^k(x))\le C \lambda^k \frac{\delta_1}4,
\quad\text{for}\quad k=0,\dots,n.
\end{align*}
Hence because $J$ is at least $C^{1+\alpha}$ for some
$\alpha\in(0,1)$ with H\"older constant $C>0$ (in fact we
can take $\alpha=1$ if $f$ is $C^2$) the usual bounded
distortion argument provides a constant $C_0>0$ such that
\begin{align*}
  \log\frac{|\det Df^n\mid F_q|}{|\det Df^n\mid F_x|}
  =
  \sum_{j=0}^{n-1} \log\frac{|\det Df(f^j(q))|}{|\det
    Df(f^j(x))|} \le
  \sum_{j=0}^{n-1} C d\big( f^j(q),f^j(x)\big)^\alpha
  \le C_0.
\end{align*}
Hence $|S_n J(q)-S_n J(x)|\le C_0$ and by the above
integration estimates we get
\begin{align*}
  \Leb\big( Q \big) \le
  \int_{\hat\gamma} C_3 e^{C_0} e^{-S_n J(x)} m\big(f^n(Q)\cap
  f^n(\beta_q)\big) \, dq \le \tilde C e^{-S_n J(x)},
\end{align*}
where $\tilde C$ is bounded by the $d_E$-dimensional area
$A_E$ of $\hat\gamma$ (which is a function of
$s<\delta_1/4$) times a uniform bound $A_F$ for the
$d_F$-dimensional area of $f^n(\beta_q)$ (which is a
function of the curvature bound $C_2$ from
Proposition~\ref{pr:bounded-curvature} and of $\delta_1$,
see Figure~\ref{fig-F-disks}) multiplied by the bounded
distortion constants, that is $\tilde C\le C_3 e^{C_0} A_E
A_F$.

This shows that we have the same kind of estimate for the
volume of a dynamical ball as in the local diffeomorphism
case, except for a different distortion constant and the
fact that the Jacobian is calculated along the $F$
direction.

\subsubsection{Weak distortion estimate}
\label{sec:weak-distort-estimat}

\textbf{If $n$ is not a hyperbolic time but
  $y\in B(x,n,\delta)$ satisfies $|S_nJ(x)-S_nJ(y)|\le n\zeta$
  %$|J(f^ky)-J(f^kx)|<\zeta$ for $k=0,\dots,n-1$
  for some pair $\delta,\zeta>0$, then
  we obtain a weak distortion property.  We can now argue similarly
  as above, using local inverse branches in the local
  diffeomorphism case (with or without singularities) to
  obtain $\Leb(B(x,n,\delta))\le C e^{n\zeta}e^{-S_nJ(x)}$;
  and $\Leb(Q)\le Ce^{n\zeta}e^{-S_nJ(x)}$ in the partially
  hyperbolic diffeomorphism case.  }
%%%%%%%%%%%%%%%%%%%%%%%%%%%%%%%%%%%%%%%%%%%%%%%%%%%%%%%

\section{Hyperbolic times and large deviations}
\label{sec:hyperb-times-large}

The statements of the main theorems and corollaries are
consequences of the following more abstract result.

\begin{theorem}
  \label{thm:hyptimes-LD}
  Let $f:M\to M$ be a local diffeomorphism outside a
  non-flat singular set $\cS$ admitting
  $\sigma\in(0,1)$ and $b,\delta>0$ such that Lebesgue
  almost every point has positive density of
  $(\sigma,\delta,b)$-hyperbolic times.  Then given
  $c\in\RR$ and a continuous function $\vfi:M\to\RR$ items
  (1)-(3) of Theorem~\ref{mthm:largedeviation} hold.
\end{theorem}

Clearly Theorem~\ref{mthm:largedeviation} follows from
Theorem~\ref{thm:tempos-hip-existem} together with
Theorem~\ref{thm:hyptimes-LD}. Moreover item (1) in the
statement of Theorem~\ref{mthm:largedeviation} is just item
(1) of \cite[Theorem 1]{Yo90} so it will not be proved here.

%%%%%%%%%%%%%%%%%%%%%%%%%%%%%%%%%%%%%%%%%%%%%%%%%%%%%%%

%%%%%%%%%%%%%%%%%%%%%%%%%%%%%%%%%%%%%%%%%%%%%%%%%%%%%%%

\subsection{Upper bound for local diffeomorphisms}
\label{sec:upper-bound-large}

Here we prove the upper bound in item 2 of
Theorem~\ref{thm:hyptimes-LD}.

Let $\vfi:M\to\RR$ be a fixed continuous function. Consider
for $n\ge1$ and some fixed $\epsilon,\delta,c>0$
\[
A_n=A_n(\delta,\epsilon)=\left\{x:
\frac1n S_n\Delta_\delta(x)\le\epsilon
\right\}
\mbox{  and  }
B_n=\left\{ x:
\frac1n S_n\vfi(x) \ge c
\right\}.
\]
Since we want to bound a limit superior from above, we can
assume without loss that $\Leb(A_n\cap B_n)>0$ in what
follows.  \textbf{We fix $\zeta>0$, set $\rho$ to be a
  positive constant function, and find
  $0<\delta_0<\delta_1/4$ and a partition $\cP$ of $M$ as in
  Subsection~\ref{sec:cover-hyperb-pre} whose diameter is
  smaller than $\de_0$ so that\footnote{Here $f$ is a local
    diffeomorphism ($\cS=\emptyset$) on a compact manifold.}
  $y\in\cP(x)\implies |J(y)-J(x)|<\zeta$}.  Then we use
Lemma~\ref{le:coverhiptimes} with $m=n$,
$E\subset U_1\subset A_n\cap B_n$ such that $U_1$ is
open\footnote{Since $S_n\vfi$ is continuous and
  $S_n\Delta_\delta$ is upper-semicontinuous.}  and
$ \Leb\big((B_n\cap A_n) \setminus E\big)<\Leb(B_n\cap
A_n)/2n.  $ Then we can find $k\ge1$ and a family
$\U_n=\cE_1\cup\dots\cup\cE_k$ of hyperbolic preballs
contained in $U_1$ satisfying
\[
\Leb\big(E\triangle \bigcup \U_n\big)
\le \left( 1-\frac\theta4 \right)^k
 < \frac1{2n} \Leb(A_n\cap B_n).
\]
Note that $\Leb\big((A_n\cap B_n)\setminus \U_n\big) \le
\Leb\big((A_n\cap B_n)\setminus E\big) + \Leb(E\setminus \U_n) < \frac1n
\Leb(A_n\cap B_n)$ and so
\begin{align}
  \label{eq:viu}
\Leb(A_n\cap B_n) < \frac{n}{n-1} \Leb(\U_n).
\end{align}
Observe also that for any element $P\in\cE_i$ there exists
$x\in M$ and a hyperbolic time $h_i$ of $f$ for $x$ such
that $P\subset B(x,h_i,\de_1)$, by construction, where
$i=1,\dots,k_n$ and $n<h_1<\dots<h_{k_n}$.  Let $\cC_n$ be
the set of all such pairs $(x,h_i)$, one for each element of
$\U_n$ and to simplify the notation we write $h_n$ for
$h_{k_n}$.

\textbf{Note that if $(x,l), (x',l')\in\cC_n$ and
  $x'\in\cP^n(x)$, then $x'\in B(x,n,2\delta_0)$ and
  $B(x,l,\delta_0)\cup B(x',l',\delta_0)\subset
  B(x,n,2\delta_0)$. Hence we may replace every pair of
  elements of $\cC_n$ in the same atom of $\cP^n$ by one of
  them, obtaining a coarser cover $U_n$ of $E$ formed by
  dynamical balls centered around a reduced family
  $\wt{\cC_n}\subset\cC_n$ of points so that, from
  Subsection~\ref{sec:weak-distort-estimat}}
\begin{align}\label{e-last}
  \Leb(A_n\cap B_n)\le\frac{n}{n-1}\Leb(U_n)
  \le
  \frac{C n e^{n\zeta}}{n-1}
  \sum \{e^{-S_n J(x)}\cdot \delta_x :
  x\in \wt{\cC_n}\}.
\end{align}
Following the arguments in the proof of \cite[Thm.1(2)]{Yo90}
\textbf{we consider the measure}
\[
  \sigma_n=Z_n^{-1}\sum \{e^{-S_n J(x)}\cdot \delta_x :
  x\in \wt{\cC_n}\}
\quad\mbox{where}\quad
Z_n=\sum\{ e^{-S_n J(x)}: x \in \wt{\cC_n}\}.
\]
\textbf{Note that by definition each atom of
$\cP^n=\bigvee_{i=0}^{n-1} f^{-i}\cP$ contains at most one point
from $\wt{\cC_n}$}. Thus using \cite[Lemma
9.9]{Wa82} we have
\begin{equation*}
  \label{e-geqlog}%  \label{eq:walters}
  H_{\sigma_n}\Big( \bigvee_{i=0}^{n-1}f^{-i}\cP \Big)
-
\int S_n J(x) \, d\sigma_n(x)
=
\log \sum\{ e^{-S_n J(x)}: x \in \wt{\cC_n}\}.
\end{equation*}
% where we write $l(x)$ for the unique integer $l$ such that
% $(x,l)\in\cC_n$. Since $S_{l(x)-n}J(f^n(x))>0$ (see
% Subsection \ref{sec:meas-dynam-balls}) and $l(x)>n$ we get
% \begin{align}
% \label{e-geqlog}
%   H_{\sigma_n}\Big( \bigvee_{i=0}^{h_n-1}f^{-i}\cP \Big)
% -
% \int S_n J \, d\sigma_n
% \ge
% \log \sum_{(x,l)\in\cC_n} e^{-S_l J(x)}.
% \end{align}
Setting $\mu_n=n^{-1}\sum_{i=0}^n f^i_*\sigma_n$ and $\mu$
a weak$^*$ accumulation point of $\mu_n$, we may modify the
initial partition $\cP$ according to
\textbf{Proposition}~\ref{pr:zeroboundary} and
Remark~\ref{rmk:partitionperturbation} so that its diameter
is smaller than $\delta_1/2$ and $\mu(\partial \cP)=0$
without loss, keeping $\cC_n$ unchanged.  As in \cite[pag.
220]{Wa82} from the above we can deduce that for every
$q\ge1$
\begin{align}
\limsup_{n\to+\infty}\frac1n\log Z_n 
&\le 
\frac1q
\limsup_{n\to+\infty} H_{\mu_n} \Big( \bigvee_{i=0}^{q-1}
f^{-i}\cP \Big) + \limsup_{n\to+\infty} \int -J\, d\mu_n
\label{eq:anteslim}
\\
&\le h_\mu(f,\cP) - \int J \, d\mu \le h_\mu(f) - \int
J \, d\mu \label{eq:nolim}
\end{align}
if $f$ is a local diffeomorphism, ensuring that $\mu$ is
$f$-invariant and that $J$ is a continuous function (in this
case $\cS=\emptyset$ and $\Delta_\delta$ plays no role, we
may take $\Delta_\delta\equiv0$ and $A_n=M$).  Observe that,
\textbf{since the points in $\wt{\cC_n}\subset\cC_n$ are contained
in $B_n$} and $\mu_n$ is a linear convex combination of
measures of the form
$n^{-1}\sum_{i=0}^{n-1}\delta_{f^i(x)}$, we get for all
$n\ge1$
\begin{align}
  \int\vfi\, \mu_n
&=
\frac1n\sum_{j=0}^{n-1} \sigma_n(\vfi\circ f^j)
=
Z_n^{-1}\sum_{x\in\wt{\cC_n}} e^{-S_n J(x)}\cdot
\frac1n\sum_{j=0}^{n-1} \vfi\big(f^j(x)\big)
\nonumber
\\
&\ge
c Z_n^{-1}\sum\{ e^{-S_n J(x)}:x\in\wt{\cC_n}\}
= c \label{eq:mediac}
\end{align}
and hence $\int\vfi\,d\mu \ge c$ also because $\vfi$ is a
continuous function.

% Note that from \eqref{eq:viu} and by
% Subsection~\ref{sec:meas-dynam-balls} we get for some
% constant $C>0$
% \begin{align}
% \Leb(B_n) 
% &\le 
% \frac{n}{n-1}\Leb(\U_n) 
% \le 
% \frac{n}{n-1}
% \sum_{(x,l)\in\cC_n}
% \Leb\Big( B(x,l,\de_1)\Big)\nonumber
% \\
% &\le
% \frac{n}{n-1}\sum_{(x,l)\in\cC_n}
% C e^{-S_l J(x)}
% % \le
% % \frac{C n}{n-1}\sum_{(x,*)\in\cC_n} e^{-S_n\xi(x)} 
% = \frac{C n}{n-1} Z_n.
% \end{align}
%because $l>n$ and $\xi\ge0$. 
Therefore we have shown that there exists $\mu\in\M_f$ such
that $\int\vfi\,d\mu\ge c$ and
\[
  \limsup_{n\to+\infty} \frac1n\log\Leb(B_n)
  \le 
\zeta+\limsup_{n\to+\infty} \frac1n\log Z_n
\le
\zeta+ h_\mu(f) - \int J\,d\mu.
\]
\textbf{Since $\zeta>0$ was arbitrary}, this completes the
proof of item 2 in the statement of
Theorem~\ref{thm:hyptimes-LD} and
Theorem~\ref{mthm:largedeviation}.

%%%%%%%%%%%%%%%%%%%%%%%%%%%%%%%%%%%%%%%%%%%%%%%%%%%%%%%

\subsection{Upper bound for partially hyperbolic diffeomorphisms}
\label{sec:upper-bound-diffeo}

Here we show that a bound similar to the one in item 2 of
Theorem~\ref{mthm:largedeviation} also holds in the case of a partially
hyperbolic non-uniformly expanding attracting set.

Let $f:M\to M$ be a diffeomorphism satisfying the conditions
of Theorem~\ref{mthm:phdiff<0}, let $\vfi:M\to\RR$ be a
continuous function, fix a real number $c$ and set
$J=\log|\det Df\mid F|$.  Observe that since we have
Lemma~\ref{le:coverpartialhyp} we may argue exactly as in
the previous subsection to arrive at an inequality just like
\eqref{e-geqlog}.

Again as in the previous subsection we consider
$\mu_n=\frac1n\sum_{i=0}^n f^i_*\sigma_n$ and $\mu$ a
weak$^*$ accumulation point of $\mu_n$. We also modify the
partition $\cP$ in such a way that the boundaries of each
atom have zero measure with respect to all measures $\mu$
and $\mu_n, n\ge1$.

The inequality~\eqref{e-geqlog} enables us to obtain
inequalities~\eqref{eq:anteslim} and~\eqref{eq:nolim}
exactly as before. Together with the volume estimates
obtained in Subsections~\ref{sec:part-hyperb-case}
\textbf{and}~\ref{sec:weak-distort-estimat} we can then arrive also
at inequality \eqref{e-last} just by using a different
distortion constant and replacing the Jacobian of $f$ by the
Jacobian of $f$ \emph{along the $F$ direction}. Hence we
obtain the upper bound given by item 2 of
Theorem~\ref{mthm:largedeviation} also in the setting of
partially hyperbolic non-uniformly expanding attracting
sets. This will be very useful to deduce
Theorem~\ref{mthm:phdiff<0} in
Subsection~\ref{sec:local-diff-case}.

%%%%%%%%%%%%%%%%%%%%%%%%%%%%%%%%%%%%%%%%%%%%%%%%%%%%%%%

\subsection{Upper bound with singular/critical set}
\label{sec:upper-bound-with}

To obtain an analogous result to (\ref{eq:nolim}) in the
limit with a transformation $f$ with non-flat singularities,
thus proving item 3 from Theorem~\ref{mthm:largedeviation}
and Theorem~\ref{thm:hyptimes-LD}, we need some extra work.
\textbf{We first use the slow recurrence condition as
  follows. Let $0\le\delta_i\le\epsilon_i$ be a sequence
  such that $\epsilon_i\searrow0$
  satisfying~\eqref{eq:SlowApprox} for all pairs
  $(\epsilon_i,\delta_ i), i\ge2$ and also
  $0<|x|<\delta_2\implies |x|^{-\beta}\le-\beta\log|x|$}.

\begin{lemma}
  \label{le:unifslowrec}
  There exists $K_0>0$ such that for any given $k\ge3$ and
  each $x\in M\setminus\cS$ satisfying
  $S_n\Delta_{\delta_2}(x)\le\epsilon_2$ and
  $S_n\Delta_{\delta_k}(x)\le\epsilon_k$, we have
  $ \sum_{i=0}^{n-1}d(x_i,\cS)^{-\beta} \le K_0n.  $
    \end{lemma}
    
    \begin{proof}
      We have
      $ \sum_{i=0}^{n-1}d(x_i,\cS)^{-\beta} \le \sum_{i\in
        B_1+B_2+B_3}d(x_i,\cS)^{-\beta}$ where
      $B_1=\{0\le i< n: d(x_i,\cS)^\beta<\delta_k\}$ and
      $B_2=\{0\le i<n: \delta_k\le d(x_i,\cS)^\beta<\delta_2\}$
      and also
      $B_3=\{0,\dots,n-1\}\setminus(B_1+B_2)=\{0\le i< n:
      d(x_i,\cS)^\beta\ge\delta_2\}$. Then, the choice of
      $x\in A_n(\delta_2,\epsilon_2)\cap
      A_n(\delta_k,\epsilon_k)$ and of $\delta_2$, ensure
      that
  \begin{align*}
    \sum_{i=0}^{n-1}\frac1{d(x_i,\cS)^\beta}
    \le
    \beta S_n\Delta_{\delta_k}(x)
    +\beta S_n\Delta_{\delta_2}(x)
    +\frac1{\delta_2^\beta}\# B_3
    \le
    \beta n\epsilon_k+\beta n\epsilon_2+n/\delta^\beta_2
    =
    K_0n
  \end{align*}
  where
  $K_0=\beta\epsilon_2+\beta\epsilon_k+\delta_2^{-\beta}\le
  2\beta\epsilon_2+\delta_2^{-\beta}$.
\end{proof}

\textbf{Then we set $\rho(x)=\exp\Delta_\delta(x)$ and use the
  non-degeneracy condition (S3) to obtain the following.}

\begin{lemma}
  \label{le:unifcontJsing}
  Given $\zeta>0$ there exists $0<\delta_0<\delta_1$ so that
  any partition $\cP$ constructed as in
  Subsection~\ref{sec:constr-an-adequate} satisfies:
  $y\in\cP^n(x)\implies |S_nJ(y)-S_nJ(x)|\le \zeta n$ for
  each $n\ge1$ and $x\in A_n(\delta_2,\epsilon_2)\cap
  A_n(\delta_k,\epsilon_k)$ for each $k\ge1$.
\end{lemma}

\begin{proof}
  The choice of $\rho$ and the construction of $\cP$ ensures
  that $y\in\cP(x)\implies d(y,x)<\rho(x)/2<d(x,\cS)/2$ for
  all $y,x\in M$. Condition (S3) for
  $y\in\cP^n(x)$ and
  $n\ge1, x\in A_n(\delta_2,\epsilon_2)\cap
  A_n(\delta_k,\epsilon_k)$ ensures
  % \footnote{We
  %   are using that $f\mid M\setminus B_\delta(\cS)$ is $C^2$
  %   and so $B$ is taken as a Lipschitz constant of
  %   $J\mid M\setminus B_{\delta_2}(\cS)$.}
  \begin{align*}
    |S_nJ(y)-S_nJ(x)|
    \le
    \sum_{i=0}^{n-1}B\frac{d(y_i,x_i)}{d(x_i,\cS)^\beta}
    \le
    B\delta_0
    \sum_{i=0}^{n-1}\frac1{d(x_i,\cS)^{\beta}}
    \le
    B\delta_0 K_0 n
  \end{align*}
  since $y_i=f^iy, x_i=f^ix$ satisfy
  $y_i\in\cP(x_i), 0\le i<n$ and by
  Lemma~\ref{le:unifslowrec}. To complete the proof we just
  have to take $0<\delta_0<\zeta(B K_0)^{-1}$.
\end{proof}

\textbf{Fixing an initial partition $\cP$ in the conditions
  of Lemma~\ref{le:unifcontJsing}, the same arguments lead
  us to \eqref{eq:viu} as in
  Subsection~\ref{sec:upper-bound-large}. Likewise, if
  $(x,l), (x',l')\in\cC_n$ and $x'\in\cP^n(x)$, then
  $x'\in B(x,n,2\delta_0)$ and
  $B(x,l,\delta_0)\cup B(x',l',\delta_0)\subset
  B(x,n,2\delta_0)$. Again, we replace every pair of
  elements of $\cC_n$ in the same atom of $\cP^n$ by one of
  them, obtaining a coarser cover $U_n$ of $E$ formed by
  dynamical balls centered around a reduced family
  $\wt{\cC_n}\subset\cC_n$ of points so that, from
  Subsection~\ref{sec:weak-distort-estimat} we again
  obtain~\eqref{e-last}.} Since the points in $\wt{\cC_n}$
are contained in $A_n\cap B_n$, \textbf{we
  reobtain}~\eqref{eq:anteslim} and~\eqref{eq:mediac}, and
also $\int \Delta_\delta \, d\mu_n \le \epsilon$ for every
$n\ge1$.
\textbf{Before we can use the same device of slightly
  perturbing $\cP$ into $\cP_{\xi,\tau}$ by
  Proposition~\ref{pr:zeroboundary}, we need to show the
  following.}

\begin{lemma}
  \label{le:Szero}
  The singular set $\cS$ has null $\mu$-measure.
\end{lemma}

\begin{proof}
Arguing by contradiction, assume that $\mu(\cS)>0$. Then
there exists $a>0$ such that $\mu\big( B(\cS,\eta) \big)\ge
a$ for all $\eta>0$. Let $\eta>0$ be chosen so that
$\mu(\partial B(\cS,\eta) \big)=0$ and
$\inf_{B(\cS,\eta)}\Delta_\delta\ge 4\epsilon/a$.

On the one hand, since $\mu$ is a weak$^*$ limit point of
$\mu_n$, there exists $n_0$ such that for $n>n_0$ we have
$\mu_n\big( B(\cS,\eta) \big)\ge a/2$. On the other hand,
since $\Delta_\delta\ge0$ we get by the choice of $\eta$
\[
\frac{4\epsilon}a \mu_n\big( B(\cS,\eta) \big)
\le
\mu_n\big(\Delta_\delta\cdot \chi_{B(\cS,\eta)}\big)
\le
\mu_n(\Delta_\delta)\le\epsilon,
\]
where $\chi_{B(\cS,\eta)}$ is the characteristic function of
$B(\cS,\eta)$, from which we deduce that $\mu_n\big(
B(\cS,\eta) \big)\le a/4$. This contradiction shows that
$\mu(\cS)=0$ and concludes the proof.  
\end{proof}

\begin{lemma}
  \label{le:intlim}
  The functions $\Delta_\delta, J$ and $\psi$ are
  $\mu$-integrable.
\end{lemma}

\begin{proof}
 Let us define the sequence of functions
\[
\Delta_\delta^k=\xi_k\circ\Delta_\delta
\mbox{  where  }
\xi_k(x)=\left\{
  \begin{array}[l]{ll}
k & \mbox{if  } |x|\ge k
\\
x & \mbox{if  } |x|<k
  \end{array}
\right.,\,\, k\ge1.
\]
For $k> k_0$ with $k_0>|\log(\delta/2)|$ and fixing
$\eta>0$, since $\Delta_\delta^k$ is continuous and
$\Delta_\delta\ge\Delta_\delta^k$ there is an integer $n_0$
such that for all $n> n_0$ we have
\[
\mu(\Delta_\delta^k) \le \mu_n(\Delta_\delta^k) + \eta
\le \mu_n(\Delta_\delta) + \eta \le \epsilon+\eta.
\]
Since this holds for all $k\ge k_0$ and
$\Delta_\delta(x)\to\infty$ when $x\to\cS$, we have proved
\[
\int_{M\setminus\cS} \Delta_\delta \, d\mu < \infty.
\]
Thus we get $\Delta_\delta\in L^1(\mu)$ since
$\mu(\cS)=0$ by Lemma~\ref{le:Szero}.

For $J$ and $\psi$, note that by conditions (S2) and
(S3) on the singular set $\cS$ it follows that there exists a
constant $\zeta>\beta$ such that on a small neighborhood $V$
of $\cS$ we have
\begin{align}\label{eq:boundnondeg}
\big|\log\|Df(x)^{-1}\|\big| + \big|\log|\det
Df(x)^{-1}|\big|
\le \zeta \big| \log d(x,\cS)\big|
\end{align}
and since $f$ is a local diffeomorphism on $M\setminus\cS$,
the $\mu$-integrability of $\Delta_\delta$ implies that of
$\psi$ and $J$. This concludes the proof of the lemma.
\end{proof}

\begin{lemma}
  \label{le:accont}
The measure $\mu$ is $f$-invariant.
\end{lemma}

\begin{proof}
  Since $\mu(\cS)=0$ by Lemma~\ref{le:Szero}, we can find a
  sequence $\eta_n\to0$ of positive numbers such that
  $\mu\big(\partial B(\cS,\eta_n)\big)=0$ for all $n\ge1$
  and $\mu \big( B(\cS,\eta_n)\big)\to0$ when $n\to\infty$.

Let us fix $\eta>0$ and a continuous function
$h:M\to\RR$. Take $n_0$ such that
\[
\mu\big( B(\cS,\eta_n)\big) \cdot \sup |h| < \eta/2
\]
for all $n>n_0$ and fix $n_1>n_0$ such that 
$
\mu\big( B(\cS,\eta_n)\big)/2
\le
\mu_n\big( B(\cS,\eta_n)\big) 
\le 2 \mu\big( B(\cS,\eta_n)\big)
$
for all $n\ge n_1$.  Then if $\tilde f$ is any continuous
extension of $f\mid M\setminus B(\cS,\eta_n)$ to $M$ (which
always exists by Tietze Extension Theorem, see e.g.
\cite{munkres}) we get
\begin{align}
  \label{eq:aproxtil}
\int \big| h\circ f - h\circ \tilde f \big| \, d\mu_n 
\le
\sup|h| \cdot \mu_n\big( B(\cS,\eta_n)\big) < \eta
\end{align}
for all $n>n_1$. Also note that (\ref{eq:aproxtil}) holds
with $\mu$ in the place of $\mu_n$.  Since $h\circ\tilde f$
is continuous there exists $n_2>n_1$ such that
$ \big| \int h\circ\tilde f \, d\mu_n - \int h\circ\tilde f
\, d\mu \big| <\eta $ for every $n>n_2$.  Hence, for $n>n_2$
we obtain that $\big|
\int h\circ\tilde f \, d\mu_n - \int  h\circ\tilde f \, d\mu
\big|$ is bounded from above by
\begin{align*}
|\mu(h\circ f) - \mu(h\circ \tilde f)|
+
|\mu(h\circ\tilde f) - \mu_n(h\circ\tilde f)|
 + |\mu_n(h\circ\tilde f) - \mu_n(h\circ f)|
\le 3\eta.
\end{align*}
Since $h$ was an arbitrary continuous function and $\eta$
was any positive number, we have shown that $f_*\mu_n\to
f_*\mu$ in the weak$^*$ topology when $n\to\infty$. This is
exactly what is needed to show that $\mu$ is $f$-invariant:
$
f_*\mu=\lim_n f_*\mu_n = \lim_n
\Big(\frac1n\sum_{j=0}^{n-1}f^j_*\sigma_n +
\frac{f^n_*\sigma_n-\sigma_n}n \Big)
=\lim_n\mu_n =\mu$ in the weak$^*$ topology,
concluding the proof.
\end{proof}

\textbf{Now we can use Proposition~\ref{pr:zeroboundary} to
  obtain a perturbation $\cP_{\xi,\tau}$ of $\cP$
  ensuring the following.}  We consider $\tilde J$ a
continuous extension of $J\chi_{M\setminus B(\cS,\xi)}$ to
$M$ with the same range (this is Tietze's Extension Theorem)
for $0<\xi<\delta$ and write
\begin{align*}
  \limsup_{n\to\infty}\mu_n(- J)
&=
\limsup_{n\to\infty}
[\mu_n\big( (-J + \tilde J)\chi_{B(\cS,\xi)}\big)+
\mu_n(-\tilde J)]
\\
&\le
2\limsup_{n\to\infty}\mu_n(\zeta\Delta_\delta) 
+ \mu(-\tilde J)
\le 2\zeta\epsilon - \mu(\tilde J)
\end{align*}
since $\tilde J$ is continuous and
$|-J + \tilde J| \chi_{B(\cS,\xi)}\le
2|J|\chi_{B(\cS,\delta)} \le 2\zeta\Delta_\delta$ by
\eqref{eq:boundnondeg}. Taking $\xi\to0$ we get
$\mu(\tilde J)\to\mu(J)$ because $J\in L^1(\mu)$ and
$H_\mu(\cP_{\xi,\tau})<\infty$, together with
\eqref{eq:anteslim} we arrive at
\[
  \limsup_{n\to+\infty}\frac1n\log Z_n \le
  h_\mu(f,\cP_{\xi,\tau}) - \int J \, d\mu +
  2\zeta\epsilon
\]
for some $\mu\in\M_f$ with $\mu(\vfi)\ge c$ and
$\Delta_\delta\in L^1(\mu)$, which is enough to prove item
(3) of Theorem~\ref{thm:hyptimes-LD} and
Theorem~\ref{mthm:largedeviation}.

%%%%%%%%%%%%%%%%%%%%%%%%%%%%%%%%%%%%%%%%%%%%%%%%%%%%%%%

\section{Strictly negative upper bound}
\label{sec:striclly-negat-upper}

Here we prove Theorem~\ref{mthm:supnegative} and
Theorem~\ref{mthm:phdiff<0}.  For a $C^1$ endomorphism $f$ it
is known \cite{Ru78} that the following inequality (also
known as \emph{Ruelle's inequality}) holds for every
$f$-invariant probability measure $\mu$
\begin{align}
  \label{eq:Ruelle}
h_{\mu}(f) \le \int \Sigma^+ \, d\mu.  
\end{align}
where $\Sigma^+$ denotes the sum of the positive
Lyapunov exponents at $\mu$-a.e. point.  In
Subsection~\ref{sec:entr-form-piec} we present a proof of
this inequality in the setting of maps which are local
diffeomorphisms away from a non-flat singular set $\cS$ with
zero Lebesgue measure, for invariant probability measures
$\mu$ such that $\log d(x,\cS)$ is $\mu$-integrable.

%\marginpar{Conferir Katok-Strelcyn!}
We note that in \cite{KS86} a similar result was proved
under more general geometric assumptions but stricter
analytic hypothesis, mostly due to the fact that in
\cite{KS86} the authors considered $M$ to be a compact
metric space admitting a finite dimensional manifold $V$ as
an open dense subset and $\cS=M\setminus V$, which demands
technical conditions on how the Riemannian metric on $V$ and
$f$ behave (including the first and second derivatives on
local charts) near $\cS$ for the proof to work. Our
conditions are similar except that we only need the
transformation $f$ to be $C^1$ but assume that $\log
d(x,\cS)$ is integrable, which is natural in our setting.

%%%%%%%%%%%%%%%%%%%%%%%%%%%%%%%%%%%%%%%%%%%%%%%%%%%%%%%

\subsection{The local diffeomorphism and partially
  hyperbolic case}
\label{sec:local-diff-case}

From Ruelle's Inequality \eqref{eq:Ruelle} and from
Subsection~\ref{sec:meas-dynam-balls} it follows that we get
a non-positive upper bound in item (2) of
Theorem~\ref{mthm:largedeviation} since $\int J\,d\mu$
equals the sum of the Lyapunov exponents of $\mu$
\cite{Os68}.  Moreover let $\mu\in\EE$ be given. Then,
\textbf{since we are assuming that each element in $\EE$ is
  a weak expanding measure}, we have
\[
\int J\,d\mu = h_\mu(f)
\le \int \Sigma^+ \, d\mu
\le \int J \, d\mu.
\]
Hence if $\mu\in\M_f$ is not in $\EE$ then the inequality
\eqref{eq:Ruelle} is strict.

To prove Theorem~\ref{mthm:supnegative} we fix a continuous
$\vfi:M\to\RR$ and % assume that there exists $\mu\in\EE$ such
% that
% \begin{align}
%   \label{eq:isolated}
% \omega_0=\inf\big\{ |\eta(\vfi)-\mu(\vfi)|: \eta\in\EE\setminus\{\mu\}
% \big\} >0
% \end{align}
% then observe that if we
replace $B_n$ in
Subsection~\ref{sec:upper-bound-large} with
\begin{align}
\label{eq:Bnovo}
B_n=\left\{ 
x\in M : 
\inf\big\{\big|
\frac1n S_n\vfi(x) - \eta(\vfi)
\big| : \eta\in\EE \big\}
> \omega
\right\}
\end{align}
for some $\omega>0$.  Then $B_n$ is an open subset of $M$
% since \eqref{eq:isolated} shows that $\EE(\vfi)$ is
% not dense in $\RR$
and we can assume without loss that $\Leb(A_n\cap B_n)>0$ in
what follows, for otherwise the limit superior
in~\eqref{e-deviationA} is smaller than any given real
number and there is nothing to prove. Hence arguing as in
Subsection~\ref{sec:upper-bound-large} we obtain a measure
$\nu\in\M_f$ satisfying $\inf\big\{|\nu(\vfi)-\eta(\vfi)|:
\eta\in\EE \big\} > \omega$, the bound of item (3) of
Theorem~\ref{mthm:largedeviation} and $\Delta_\delta\in
L^1(\nu)$ with $\nu(\Delta_\delta)\le\epsilon$.

If $f$ is a local diffeomorphism, i.e.  $\cS=\emptyset$,
then we can use the bound given by item (2) of
Theorem~\ref{mthm:largedeviation} and it is enough to show
that $h_\nu(f)-\nu(J)$ is strictly negative.  But we
cannot have $h_\nu(f) - \nu(J)=0$ since by construction
$\nu$ is not in $\EE$, thus $h_\nu(f)-\nu(J)<0$,
completing the proof of Theorem~\ref{mthm:supnegative} in
the case of a local diffeomorphism.

%%%%%%%%%%%%%%%%%%%%%%%%%%%%%%%%%%%%%%%%%%%%%%%%%%%%%%%

For a partially hyperbolic non-uniformly expanding
attracting set we obtain a negative upper bound following
the same reasoning as above since we can use the same bound
from item (2) of Theorem~\ref{mthm:largedeviation}, as shown
in Subsection~\ref{sec:upper-bound-diffeo}, and we can also
apply Ruelle's Inequality. This completes the proof of
Theorem~\ref{mthm:phdiff<0}.

%%%%%%%%%%%%%%%%%%%%%%%%%%%%%%%%%%%%%%%%%%%%%%%%%%%%%%%

\subsection{The case with singular/critical set}
\label{sec:case-with-sing}

In the case $\cS\neq\emptyset$ we now show that the upper
bound in item (3) of Theorem~\ref{mthm:largedeviation} must
be strictly negative for some values of
$\eta,\epsilon,\delta>0$ and for some $\nu\in\M_f$. For that
we argue by contradiction and take \textbf{decreasing
  sequences $\delta_k\le\epsilon_k, k\ge2$ such that
  $\epsilon_k\searrow0$, each pair $(\delta_k,\epsilon_k)$
  satisfies~\eqref{eq:SlowApprox} and assume that} the corresponding
measures $\nu_k$ obtained according to the proof of
Theorem~\ref{mthm:largedeviation}, with $B_n$ as in
\eqref{eq:Bnovo} and
$ A_n^k=\cap_{i=2}^k A_n(\delta_i,\epsilon_i)=\{ x\in M:
S_n\Delta_{\delta_i}\le n\epsilon_i, i=2,\dots,k\} $ in the
place of $A_n$, for each $k\ge2$, also satisfy
\begin{itemize}
\item $\nu_k\in\M_f$, $\Delta_{\delta_i}\in L^1(\nu_k)$ and
  $\nu_k(\Delta_{\delta_i})\le\epsilon_i$ for $i=1,\dots,k$;
\item $\limsup_{n\to\infty}\frac1n\log\Leb(A_n^k\cap B_n)\le
  h_{\nu_k}(f,\cP)-\int J\,d\nu_k + 2\zeta\epsilon_k$;
\item $h_{\nu_k}(f,\cP)-\int J\,d\nu_k +
  2\zeta\epsilon_k\ge0$; and
\item
  $\inf\big\{|\nu_k(\vfi)-\eta(\vfi)|: \eta\in\EE \big\} >
  \omega$.
\end{itemize}
Above, $\cP=\cP_{\xi,\tau}$ is a partition obtained using
\textbf{Proposition}~\ref{pr:zeroboundary} with the
sequence\footnote{Here the independence of $K_0$ from
  $k\ge2$ in Lemma~\ref{le:unifslowrec} is crucial, allowing
  the size of the partition $\cP$ not to shrink with $k$.}
$\mu_k=\nu_k$ and $\mu$ some weak$^*$ accumulation point of
the $\nu_k$.  Thus, on the one hand, we have for any fixed
$N\ge1$
\[
h_{\nu_k}(f,\cP)=\inf_{j\ge1}\frac1j
H_{\nu_k}\left(\bigvee_{i=0}^{j-1} f^{-i}\cP\right)
\le
\frac1N H_{\nu_k}\left(\bigvee_{i=0}^{N-1} f^{-i}\cP\right)
\]
and since $\mu(\partial\cP)=0$ we get
\[
\limsup_{k\to\infty}h_{\nu_k}(f,\cP)
\le
\frac1N H_{\mu}\left(\bigvee_{i=0}^{N-1} f^{-i}\cP\right).
\]
But $N\ge1$ was arbitrarily fixed, so
\[
\limsup_{k\to\infty}h_{\nu_k}(f,\cP)
\le\inf_{N\ge1}
\frac1N H_{\mu}\left(\bigvee_{i=0}^{N-1} f^{-i}\cP\right)
=
h_\mu(f,\cP).
\]
On the other hand, choosing $J_i$ to be a continuous
extension of $J\chi_{M\setminus B(\cS,\delta_i)}$ to $M$ with the
same range, $i\ge1$, we have
\begin{align*}
  \limsup_{k\to\infty}\nu_k(- J)
&=
\limsup_{k\to\infty}
[\nu_k\big( (-J+J_i)\chi_{B(\cS,\delta_i)}\big)+
\nu_k(-J_i)]
\\
&\le
2\limsup_{k\to\infty}\nu_k(\zeta\Delta_{\delta_i}) 
+ \mu(-J_i)
\le 2\zeta\epsilon_i - \mu(J_i)
\end{align*}
since $J_i$ is continuous and
$|-J+J_i|\chi_{B(\cS,\delta_i)}\le
2|J|\chi_{B(\cS,\delta_i)} \le 2\zeta\Delta_{\delta_i}$
by definition of $\Delta_{\delta_i}$ and
by~\eqref{eq:boundnondeg}. Similar arguments to the ones
proving Lemmas~\ref{le:Szero},~\ref{le:intlim}
and~\ref{le:accont} show that $J,\psi,\Delta_\delta$ are
$\mu$-integrable and that $\mu$ is $f$-invariant.  Because
$i\ge1$ can be arbitrarily chosen above and both
$\epsilon_i\to0$ and $\mu(J_i)\to\mu(J)$, we conclude
that $\limsup_{k\to\infty}\nu_k(- J) \le - \mu(J)$.
Hence we deduce
\[
0\le\limsup_{k\to\infty}
\left(h_{\nu_k}(f,\cP)+\nu_k(-J)+2\zeta\epsilon_k\right)
\le
h_\mu(f,\cP)-\mu(J) \le h_\mu(f)-\mu(J)
\]
and also that $\inf\big\{|\mu(\vfi)-\eta(\vfi)|: \eta\in\EE
\big\}\ge\omega>0$ by construction. By Ruelle's Inequality
we also get $h_\mu(f)-\mu(J)\le0$, which yields a
contradiction since this means $\mu\in\EE$.  This
contradiction shows that for some $k\ge2$
\[
h_{\nu_k}(f,\cP)-\int J\,d\nu_k + 2\zeta\epsilon_k<0.
\]
This proves Theorem~\ref{mthm:supnegative}, except for the
Ruelle Inequality for maps with non-flat singularities,
which is the content of the next subsection.

%%%%%%%%%%%%%%%%%%%%%%%%%%%%%%%%%%%%%%%%%%%%%%%%%%%%%%%

\subsection{Ruelle's Inequality for maps with non-flat
singularities}
\label{sec:entr-form-piec}

\begin{theorem}
  \label{thm:Ruelle}
  Let $f:M\setminus\cS\to M$ be a $C^1$ local diffeomorphism
  away from a non-flat singular set $\cS$ and $\mu$ a
  $f$-invariant probability measure such that $|\log
  d(x,\cS)|$ is $\mu$-integrable.  Then
\[
h_\mu(f)\le \int \Sigma^+ \,d\mu,
\]
where $\Sigma^+$ denotes the sum of the positive
Lyapunov exponents at a regular point, counting
multiplicities.
\end{theorem}

Observe that the $\mu$-integrability of $|\log d(x,\cS)|$
implies the $\mu$-integrability of $\log^+\|Df\|$, where
$\log^+x=\max\{0,\log x\}$, and thus the Lyapunov exponents
of $f$ are well defined $\mu$-almost everywhere by
Oseledets Theorem \cite{Os68}.  The proof we present here
follows Ma\~n\'e \cite[Chap. IV]{Man87} closely.

  We start by taking the $M$ as a compact submanifold of
  $\RR^N$ with the usual Euclidean norm and induced
  Riemannian structure, and considering $W_0$ an open
  \emph{normal tubular neighborhood} of $M$ in $\RR^N$, that
  is, there exists $\Phi:W_0\to W, (x,u)\mapsto x+u$ a
  ($C^\infty$) diffeomorphism from a neighborhood $W_0$ of the
  zero section of the normal bundle $TM^\perp$ of $M$ to
  $W$. Let also $\pi:W\to M$ be the associated
  projection: $\pi(w)$ is the closest point to $w$ in $M$
  for $w\in W$, so that the line through the pair of
  points $w,\pi(w)$ is normal to $M$ at $\pi(w)$, see e.g.
  \cite{hirsch1976} or \cite{guillemin-pollack1974}.  Now we
  define for $\rho\in(0,1)$
\[
F_0:W_0\setminus(T_{\cS}M)\to W_0,\quad (x,u)\mapsto (f(x),\rho\cdot u)
\]
and also
\[
F:W\setminus\Phi(T_{\cS}M)\to W, \quad w\mapsto (\Phi
\circ F_0 \circ \Phi^{-1})(w).
\]
Then clearly $F$ is a local diffeomorphism outside
$\Phi(T_{\cS}M)$, $\overline{F(W)}\subset W$ and
$M=\cap_{n\ge0} F^n(W)$.

% Define $\tilde\cP=\{W\setminus\Phi(T_{\cS}M), W\}$ a
% partition of $W$, 
For each $n\ge1$ consider the partition
of $\RR^N$ into dyadic cubes
\[
\cP_n=\left\{\prod_{i=1}^N
  \Big[\frac{a_i}{2^n},\frac{a_i+1}{2^n}\Big)
: a_i\in\ZZ, i=1,\dots, N\right\}.
\]
%and set $\cP_n=\tilde\cP\vee\tilde\cP_n$. 
Up to a slight translation of the partitions $\cP_n$ we can
assume that the probability measure $\mu$ on $M$ satisfies
$\mu(M\cap\partial\cP)=0$, where
$\partial\cP=\cup_{n\ge1}\partial\cP_n\cup\cS$. For $x\in
M\setminus\partial\cP$ we define 
\[
v_n(x)=v^F_n(x)=\#\{P\in\cP_n:F(\cP_n(x))\cap
P\ne\emptyset\} 
\]
and 
\[
v(x)=v^F(x)=\limsup_{n\to\infty}
v_n(x)
\]
where $\cP_n(x)$ denotes the atom of the partition $\cP_n$
containing $x$.

\begin{lemma}
  \label{le:Mane12.1}
Let $Q=[-1,1]^N$ and $x\in
M\setminus\partial\cP$. Then
\[
v(x)\le\sup_{z\in\RR^n}\#\{P\in\cP_1:
\big(z+Dg(x)Q\big)\cap P\ne\emptyset\}
\]
\end{lemma}

\begin{proof}
  For $x\in M\setminus\partial\cP$ and $n\ge1$ define
  $\vfi_n(y)=x+y/n$ on $\RR^N$ and $W_n=\vfi_n^{-1}(W)$. Let
  $F_n:W_n\to F_n(W_n)\subset W_n$ be such that
\[
\begin{array}[c]{rcl}
W_n & \stackrel{F_n}{\longrightarrow} & W_n
\\
\vfi_n \downarrow  & & \downarrow\vfi_n
\\
W & \stackrel{F}{\longrightarrow} & W
\end{array}
\]
commutes. We have $F(w)=F(x)+DF(x)(w-x)+p_x(w)$ where
$p_x:W\setminus\Phi(T_{\cS}M)\to\RR^N$ is $C^1$ and
$\lim_{w\to x}\|p_x(w)\|/\|w-x\|=0$, where $\|\cdot\|$ is
the Euclidean norm on $\RR^N$. Then we write
$F_n(y)=DF(x)(y)+q_n^x(y)+\alpha_n(x)$ where
\begin{align}\label{eq:qnx}
\alpha_n(x)=n\cdot F(x)-x
\quad\mbox{and}\quad
q_n^x(y)=n\cdot p_x\big(y/n+x).
\end{align}
Note that for $x\in M\setminus\partial\cP$ we have
$q_n^x\to0$ uniformly on compacta. Indeed if $\|y\|<r$ for
some $r>0$ there is, for each given $\delta>0$, a
$n_0\in\NN$ such that $\| y/n \|<\delta,\forall n\ge n_0$
and then, by definition of $p_x$, for all $\epsilon>0$ there
is $n_1\in\NN$ so that $\forall n\ge n_1,
\|p_x(y/n+x)\|<\epsilon\|y/n\|$ which is the same as
$\|n\cdot p_x(y/n+x)\|<\epsilon r$, or
$\|q_n^x(y)\|<\epsilon r$ for all sufficient large $n$.

Commutativity of the diagram implies $$
F(\cP_n(x))\cap
P\ne\emptyset\Leftrightarrow
F_n(\varphi_n^{-1}(\cP_n(x)))\cap\varphi_n^{-1}(P)\ne\emptyset.
$$
But $\varphi_n^{-1}(P)$ is an element of $\cP_1$
translated by some vector $y_0\in\RR^N$. Moreover
$\varphi_n^{-1}(\cP_n(x))\subset Q$ and so 
$
v_n(x)\le\#\{P\in\cP_1:F_n(Q)\cap(P+y_0)\ne\emptyset\}.
$
Because $\alpha_n$ depends on $x$ only 
\begin{align}
v_n(x) 
&\le
\#\left\{P\in\cP_1:\Big(n\cdot DF(x)(\frac1{n}Q)+q_n^x(Q)+
\alpha_n(x)-y_0\Big)\cap P\ne\emptyset
\right\}
\nonumber
\\
&\le
\sup_{z\in\RR^N}\#
\left\{P\in\cP_1:\big(DF(x)Q+q_n^x(Q)+ z\big)\cap P
\ne\emptyset\right\}\label{eq:boundvn}
\end{align}
Since $q_n^x\to0$ on compact subsets we get
\[
\limsup_{n\to\infty} v_n(x) \le
\sup_{z\in\RR^N}\#
\left\{P\in\cP_1:\big(DF(x)Q+ z\big)\cap P
\ne\emptyset\right\}
\]
concluding the proof of the lemma.
\end{proof}

For the arguments which use the convergence properties of
the sequence $\log v_n$ we need the following result.

\begin{lemma}
  \label{le:dominated}
  There exists a $\mu$-integrable function $g$ such that
  $0\le\log v_n\le g$ for $\mu$-almost every point in $M$
  and for all $n\ge1$.
\end{lemma}

\begin{proof} 
  Fix $n\ge1$ and consider $x\in M\setminus\partial\cP$.  On
  the one hand since $\cP_n$ is a partition we must have
  $v_n(x)\ge1$. On the other hand, by the bound
  \eqref{eq:boundvn} since the size of the edge of the
  cubes of $\cP_1$ is $1/2$ in $\RR^N$ we get
  \begin{align}
    v_n(x)
    &\le
    \Big( 2\big(\diam DF(x)(Q) + \diam q_n^x(Q) \big)
    \Big)^N \label{eq:vndiam}
    \\
    \diam DF(x)(Q)
    &\le 
    2\sqrt{N}\cdot\|DF(x)\| \nonumber
    \\
    &\le
    2\sqrt{N}\max\{\|Df(x)\|,\| DF\mid \big(T_x M)^\perp
    \|\}. \label{eq:diamDFQ}
  \end{align}
  Note that for $x$ far away from $\cS$ we always get
  bounded expressions above since $F$ is a local
  diffeomorphism outside of $\Phi(T_\cS M)$.  To bound
  $\diam q_n^x(Q)$ we use \eqref{eq:qnx} and consider two
  cases.

First assume that $d(x,\cS)\ge2/n$ and take $y\in Q$. Then
for some $\theta\in[0,1]$
\begin{align*}
  q_n^x(y)
  &=
  n\cdot p_x(y/n+x)
  = n\cdot\Big( F(x+y/n) - F(x) - DF(x)(y/n) \Big)
%   \\
%   &=
%   n\cdot\Big( DF(x+\theta \cdot y/n )(y/n) - DF(x)(y/n)
%   \Big)
  \\
  &= DF(x+\theta \cdot y/n )(y) - DF(x)(y)
\end{align*}
so we get by condition (S1) on $\cS$ 
\begin{align}
\|q_n^x(y)\| 
&\le 
\sqrt{N}\cdot\big( \|DF(x)\| + \|DF(x+\theta \cdot y/n )\|
\big)\nonumber
\\
&\le
B\sqrt{N}\Big( d(x,\cS)^{-\beta}
+\big(d(x,\cS)-1/n\big)^{-\beta} \Big)\nonumber
\\
&\le
B\sqrt{N} \cdot d(x,\cS)^{-\beta} \cdot (1+2^\beta)\label{eq:qnx1}
\end{align}
since $1-1/(nd(x,\cS))\ge 1/2$ and $\| DF\mid \big(T_x
M)^\perp \| \le\rho<1\ll d(x,\cS)^{-\beta}$ for $x$ close to
$\cS$, because $\beta>0$.

Now assume that $d(x,\cS)<2/n$. Then we bound as follows
\begin{align}
  \|q_n^x(y)\| 
  &\le 
  n\cdot\| F(x+y/n) - F(x) \| + \| DF(x) \|\cdot\|y\|\nonumber
  \\
  &\le
  n\cdot\diam W + B\sqrt{N}\cdot d(x,\cS)^{-\beta}\label{eq:qnx2}
\end{align}
Hence putting \eqref{eq:vndiam}, \eqref{eq:diamDFQ},
\eqref{eq:qnx1} and \eqref{eq:qnx2} together we see that
there exists a constant $\tilde C>0$ such that
\[
\log v_n(x) \le \left\{
  \begin{array}[l]{ll}
    N \log \big( \tilde C d(x,\cS)^{-\beta}\big) 
    & \mbox{  if  } d(x,\cS)\ge2/n,
    \\
    N \log \big( \tilde C d(x,\cS)^{-\beta}  + 2n\cdot\diam W \big)
    & \mbox{  if  } d(x,\cS)<2/n.
  \end{array}
\right.
\]
But $d(x,\cS)^{-\beta}>0$ and we may assume
without loss that $2n\cdot\diam W\ge2$, so
\[
\log \big( \tilde C d(x,\cS)^{-\beta}  + 2n\cdot\diam W \big)
\le  
\log \big( \tilde C d(x,\cS)^{-\beta} \big)
+ \log \big( 2n\cdot\diam W\big)
\]
and if $d(x,\cS)<2/n$ we also get
\begin{align*}
  \log d(x,\cS)^{-\beta}
  &=
  -\beta \log d(x,\cS) \ge -\beta \log (2/n) = \beta\log (n/2)
  \\
  &=
  \beta\log\big(  2n\cdot\diam W \big) -\beta\log(4 \diam W)
  \quad\mbox{or}
  \\
  \log \big( 2n\cdot\diam W\big)
  &\le
  \log(4 \diam W) - \log d(x,\cS)
\end{align*}
Hence in all cases we arrive at
\[
\log v_n(x) \le N\log\Big( C d(x,\cS)^{-\beta} + D\Big)
\]
for some positive constants $C$ and $D$. This concludes the
proof.
\end{proof}

\begin{lemma}
  \label{le:Mane12.2}
The following bound on the entropy holds
\[
h_\mu(f,\cP_n\cap M )=
h_\mu\big( F\mid M , \cP_n\cap M \big)
\le
\int_M \log v^F_n \, d\mu.
\]
\end{lemma}

\begin{proof}
  This is \cite[Lemma 12.2]{Man87} without change.
\end{proof}

\begin{corollary}
  \label{cor:Mane12.2}
$h_\mu(f)=h_\mu(F\mid M) \le \int_M \log v^F \, d\mu$.
\end{corollary}

\begin{proof}
  Since $\bigvee_{n\ge1}(\cP_n\cap M)$ is the Borel
  $\sigma$-algebra $\mu\bmod0$ we get
\[
h_\mu(F\mid M)=
\lim_{n\to\infty} h_\mu\big( F\mid M , \cP_n\cap M \big)
\le
\limsup_{n\to\infty} \int_M \log v^F_n \, d\mu.
\]
By Lemma~\ref{le:dominated} we can use the Dominated
Convergence Theorem to obtain
\[
\limsup_{n\to\infty} \int_M \log v^F_n \, d\mu
\le
\int_M \limsup_{n\to\infty} \log v^F_n \, d\mu
=
\int_M \log v^F \, d\mu
\]
since $\log$ is monotonous increasing. This concludes the
proof.
\end{proof}

In what follows write $v^n(x)=v^{F^n}(x)$ for the analogous
to $v^F(x)$ with $F^n$ in the place of $F$.

\begin{lemma}
  \label{le:Mane12.3}
We have
\[
h_\mu(f)=h_\mu(F\mid M) \le
\int\limsup_{n\to\infty}\frac1n\log v^n(x)\,d\mu(x).
\]
\end{lemma}

\begin{proof}
  Using \cite[Thm. 4.13]{Wa82} and Corollary~\ref{cor:Mane12.2}
  we get for all $n\ge1$
  \begin{align}\label{e-hFn}
    h_\mu(F\mid M)=\frac1n h_\mu(F^n\mid M) \le
    \int\frac1n\log v^n(x)\,d\mu(x).
  \end{align}
Consider the sequence $g_n(x)=n^{-1}\log v^n(x)$ and observe
that by Lemma~\ref{le:Mane12.1} and by~\eqref{eq:diamDFQ}
\begin{align}
  g_n(x) & \le \frac1n\log\big( 2\diam(DF^n(x)Q) \big)^N
  \nonumber
  \\
  &\le\frac{N}n\log(2\sqrt{N}) +
  \frac{N}n\log\|DF^n(x)\|=G_n(x).\label{e-gn0}
\end{align}
Again by~\eqref{eq:diamDFQ} and by definition of $F$ since
$x\in M$ we get $\log\|DF(x)\|\le \log^+\|Df(x)\|$.  Hence
by the $f$-invariance of $\mu$ and the Sub-additive Ergodic
Theorem \cite[Thm. 10.1]{Wa82}, the sequence $G_n(x)$ tends
to a finite limit $G(x)$ for $\mu$-a.e. $x$ when
$n\to\infty$.

Now by~\eqref{e-gn0} and by Fatou's Lemma~\cite[Thm.
0.9]{Wa82}
\begin{align}\label{e-Fatou}
  \int \liminf_{n\to\infty} (G_n-g_n) \, d\mu
  \le
  \liminf_{n\to\infty} \int (G_n-g_n) \, d\mu.
\end{align}
On the one hand since $\lim_{n\to\infty}
G_n(x)$ exists $\mu$-a.e.
\begin{align}\label{e-liminf}
  \int \liminf_{n\to\infty} (G_n-g_n) \, d\mu
  = \int (G-\limsup_{n\to\infty} g_n)\, d\mu
\end{align}
and, on the other hand, since $\lim_{n\to\infty}\int
G_n(x)\,d\mu$ exists $\mu$-a.e. we also get
\begin{align}
  \label{e-limsup}
  \liminf_{n\to\infty} \int (G_n-g_n) \, d\mu
  =
  \int G\,d\mu - \limsup_{n\to\infty} \int g_n \, d\mu.
\end{align}
Altogether \eqref{e-Fatou}, \eqref{e-liminf} and
\eqref{e-limsup} imply 
\begin{align*}
\limsup_{n\to\infty}\int\frac1n\log v^n(x)\,d\mu(x)
 \le
 \int\limsup_{n\to\infty}\frac1n\log v^n(x)\,d\mu(x)
 \end{align*}
 which together with~\eqref{e-hFn} conclude the proof of
 the Lemma.
\end{proof}
To finish we need to relate $\limsup_{n\to\infty}\frac1n\log
v^n(x)$ with the sum of the positive Lyapunov exponents at
$x$. This is done just as in \cite[Chap. IV, Sec. 12]{Man87}
where it is proved that
\[
\limsup_{n\to\infty}\frac1n\log
v^n(x) \le \Sigma^+(x)
\]
for $\mu$-almost all $x\in M$. This together with
Lemma~\ref{le:Mane12.3} implies Ruelle's Inequality. The
proof of Theorem~\ref{thm:Ruelle} is complete.

\bibliographystyle{abbrv}
%\bibliography{../../Trabalho/bibliobase/bibliography}

\section{List of changes with respect to the published
  version}
\label{sec:lista-changes-with}

Compared to Journal of Statistical Physics,
Vol. 125(2):415--457, 2006 (DOI:
10.1007/s10955-006-9183-y), in this version we have written
in boldface the following changes: we have
\begin{enumerate}
\item added the assumption that all equilibrium states
  $\mu\in\EE$ are weak expanding in the statements of
  Theorems~\ref{mthm:supnegative} and~\ref{mthm:phdiff<0}
  and their Corollaries~\ref{mcor:escaperate}
  and~\ref{mcor:escaperatediffeo}.
\item corrected comments on the examples, since some of them
  do not satisfy the extra assumption introduced above.
\item added a reference to a related paper on exponentially
  slow approximation to the singular set of Lorenz-like maps
  which was published at a later date.
\item promoted Lemma~\ref{pr:zeroboundary} to
  Proposition~\ref{pr:zeroboundary}; strengthen its
  statement and provide full proof of the new statement.
\item added Subsection~\ref{sec:weak-distort-estimat} to
  state a weak bounded distortion estimate that will be
  used.
\item completed the argument providing the
  bounds~\eqref{e-geqlog} and~\eqref{e-last} in the local
  diffeomorphism and partial hyperbolic cases.
\item adapted the argument in the singular/critical case in
  Subsection~\ref{sec:upper-bound-with} introducing new
  Lemmas~\ref{le:unifslowrec} and~\ref{le:unifcontJsing}
  crucially using the new Proposition~\ref{pr:zeroboundary}
  to obtain a complete proof of
  Theorem~\ref{mthm:largedeviation}.
\end{enumerate}
Finally, we updated the affiliation and emails of the first author.

\end{document}